\documentclass{amsart}
\usepackage{algorithm}
\usepackage{algorithmic}
\usepackage{graphics}
\usepackage{color}
\usepackage{amsmath}
\usepackage{amsthm}
\usepackage{amsfonts}
\newtheorem{theorem}{Theorem}
\newtheorem{lemma}[theorem]{Lemma}
\newtheorem{proposition}[theorem]{Proposition}
\newtheorem{corollary}[theorem]{Corollary}
\newcommand{\baseRing}[1]{\ensuremath{\mathbb{#1}}}
\newcommand{\N}{\baseRing{N}}

\newcommand{\R}{\baseRing{R}}

\newcommand{\Q}{\baseRing{Q}}

\newcommand\gameenv[1]{$(-#1)$-strategy profile}
\newcommand\V[1]{$({\rm\bf V}#1)$}
\newcommand\E[1]{$({\rm\bf E}#1)$}
\newcommand\startsystem[1]{$({\rm\bf S}#1)$}
\newcommand\sign{\hbox{\rm sign}}
\newcommand\Gambit{{\sf Gambit}}
\newcommand\PHC[1]{{\sf PHC}}
\newcommand\Singular{{\sf Singular}}
\newcommand\JeronimoPerrucciSabia{Jeronimo, Perrucci, and Sabia}
\begin{document}
\title{Finding All Nash Equilibria of a Finite Game Using Polynomial Algebra}
\author{Ruchira S. Datta}
\email{datta@math.berkeley.edu}
\urladdr{http://math.berkeley.edu/\~datta}
\keywords{Nash equilibrium, normal form game, algebraic variety}
\date{\today}
\begin{abstract}
The set of Nash equilibria of a finite game is the set of nonnegative solutions
to a system of polynomial equations.  In this survey article we describe how to
construct certain special games and explain how to find all the complex roots
of the corresponding polynomial systems, including all the Nash equilibria.  We
then explain how to find all the complex roots of the polynomial systems for
arbitrary generic games, by polyhedral homotopy continuation starting from the
solutions to the specially constructed games.  We describe the use of Gr\"obner
bases to solve these polynomial systems and to learn geometric information
about how the solution set varies with the payoff functions.  Finally, we
review the use of the \Gambit\ software package to find all Nash equilibria of
a finite game.
\end{abstract}
\maketitle
\section{Introduction}
\noindent
The set of Nash equilibria of a finite game is the set of nonnegative solutions
to a system of polynomial equations.  In this article we introduce this point
of view and survey the computational methods for finding Nash equilibria using
polynomial algebra which are available to game theorists today.  We give
examples which we hope will enable game theorists to find all Nash equilibria
of games of relatively larger formats.  When we know only a subset of the Nash
equilibria of a game, we have a rather shaky foundation for making predictions
about what could happen in repeated trials of the game, since play may converge
toward a Nash equilibrium which is not in the subset.  Knowing all the Nash
equilibria of a game will help us to make more principled predictions.

\section{The System of Polynomial Equations}
\noindent
A \emph{monomial} in $n$ variables $x_1,\ldots,x_n$ is an expression of the
form $x_1^{d_1}\ldots x_n^{d_n}$ for some nonnegative integers
$d_1,\ldots,d_n$.  It is \emph{squarefree} if $d_i\leq 1$ for all $i$.  The
\emph{degree in $x_i$} of this monomial is $d_i$, and its \emph{total
degree} is $\sum_{i=1}^n d_i$.  The product of two monomials
$x_1^{d_1}\ldots x_n^{d_n}$ and $x_1^{e_1}\ldots x_n^{e_n}$ is
$x_1^{d_1+e_1}\ldots x_n^{d_n+e_n}$.

A \emph{polynomial} in $n$ variables $x_1,\ldots,x_n$ with coefficients in
a field $K$ is a finite sum of \emph{terms}.  Each term is of the form $cm$
for some $c\in K$ and some monomial $m$ in $x_1,\ldots,x_n$.  The set of
polynomials in $n$ variables is a vector space whose basis is the set of
monomials in $n$ variables.  Since we can also multiply monomials together,
we can extend this multiplication to define a product on the space of
polynomials, using commutativity, associativity, and distributivity.  This
makes the set of polynomials in $n$ variables with coefficients in $K$ into a
\emph{commutative ring}, and the study of such objects is the subject of
\emph{commutative algebra}.  A \emph{field} is a commutative ring in which
every nonzero element has a multiplicative inverse.

We can \emph{evaluate} a polynomial
$f=\sum c_{d_1\ldots d_n}x_1^{d_1}\cdots x_n^{d_n}$ at a point $(a_1,\ldots,a_n)\in K^n$ by
\emph{substituting} $a_1,\ldots,a_n$ for the variables $x_1,\ldots,x_n$
respectively, obtaining an expression $\sum c_{d_1\ldots
d_n}a_1^{d_1}\cdots a_n^{d_n}$ and carrying out all the multiplications and
additions in the field $K$.  We denote the resulting element of $K$ as
$f(a_1,\ldots,a_n)$.  

A \emph{polynomial equation} is an expression $f=g$
for some polynomials $f$ and $g$.  A point $a\in K^n$ satisfies this
equation if $f(a)=g(a)$.  Since this is equivalent to $f(a)-g(a)=0$, we can
always write a polynomial equation as $p=0$ for some polynomial $p$.  A
point satisfying this equation is called a \emph{root} of $p$.  A
polynomial in one variable with coefficients in $K$ need not have any roots
in $K$.  But there is always a field containing $K$, called the
\emph{algebraic closure $\bar K$ of $K$}, such that every nonconstant
univariate polynomial with coefficients in $\bar K$ has a root in $\bar K$.
The field of complex numbers is the algebraic closure of the field of real
numbers.  

The \emph{multiplicity} of a root of a polynomial denotes how many linear
factors of that polynomial vanish at that root.  (Every polynomial factorizes
into linear factors over the algebraic closure of its coefficient field.)  For
example, the multiplicity of the root $1$ of the polynomial $(x-1)^2$ is $2$.
A polynomial equation of degree $d$ in one variable has $d$ complex roots
(counted with multiplicity); this is the \emph{Fundamental Theorem of Algebra}.  

A \emph{system of polynomial equations in $n$ variables over $K$} is a
finite set of polynomial equations in $n$ variables over $K$, and a point
in $K^n$ satisfies, or is a root of, this system if it satisfies all the
constituent equations.  The study of solution sets of polynomial systems
is called \emph{algebraic geometry}.  Recent years have seen a renaissance
in \emph{computational algebraic geometry}; the interested reader is
referred for example to \cite{CLO1}.  The subsequent volume \cite{CLO2}
contains more information, particularly about finding the roots of polynomial
systems.  The book \cite{sturmfels} surveys techniques for solving polynomial
systems and includes a chapter on finding Nash equilibria.
\cite{DickensteinEmiris} is a recent summary of the state of the art of
solving polynomial systems.  \cite{SommeseWampler} explains numerical algebraic
geometry and is written for a general technical audience (rather than a
mathematical one).

We now fix the game-theoretic notation we shall use.  The concepts we describe
here can be found in a standard game theory text such as
\cite{OsborneRubinstein}.  A normal form game with a
finite number of players, each with a finite number of pure strategies, is
specified as follows.  The set of players is denoted as $I=\{1,\ldots,N\}$.
Associated to the players are finite disjoint sets of pure strategies
$S_1,\ldots,S_N$.  We require that $|S_i|\geq 2$ for each $i$.
Write $S=\prod_{i\in I}S_i$ for the set of pure strategy
profiles.  For each $i$ let $d_i=|S_i|-1$, and write the set $S_i$ as
$\{s_{i0},\ldots,s_{id_i}\}$.  So a pure strategy profile $s\in S$ can be
written as $s=(s_1,\ldots,s_N)$, where for each $i$, we have $s_i=s_{ij}$ for
some $j$ with $0\leq j\leq d_i$.  Let $D=\sum_{i\in I}\left(|S_i|-1\right)$.
So $\sum_{i\in I}|S_i|=D+N$.
The set $\Sigma_i$ of mixed strategies of
player $i$ is the set of all functions $\sigma_i\colon S_i\to [0,1]$ with
$\sum_{j=0}^{d_i}\sigma_i(s_{ij})=1$.  Write $\Sigma=\prod_{i\in I}\Sigma_i$
for the set of strategy profiles. Write $\Sigma_{-i}=\prod_{j\in I-\{i\}}
\Sigma_j$.  We will call an element of $\Sigma_{-i}$ a \gameenv{i}.  Write
$S_{-i}=\prod_{j\in I-\{i\}}S_j$.  We will call an element of $S_{-i}$ a pure
\gameenv{i}.  Write
$\sigma_{-i}$ for the image of $\sigma\in\Sigma$ under the projection
$\pi_{-i}$ onto $\Sigma_{-i}$.  Write $\sigma_{ij}=\sigma_i(s_{ij})$ for
$j=0,\ldots,d_i$.  A \emph{game format} is a specification of a number of
players $N$ and the number of pure strategies each player has, specified by
the numbers $d_1,\ldots,d_N$.  Without loss of generality we can require that 
$d_1\geq d_2\geq\cdots\geq d_N$.

The game is specified by describing the
payoff function $u_i\colon S\to\R$ for each player.  The $i$th player's
expected payoff from a strategy profile $\sigma$ is given by multilinearity as
$$u_i(\sigma)=\sum_{s\in S} u_i(s)\sigma_1(s_1)\cdots\sigma_N(s_N).$$
By abuse of notation, write $u_i(\sigma_i,\sigma_{-i})$ for the $i$th
player's expected payoff from the strategy $\sigma$ whose $i$th component
is $\sigma_i$ and whose other components are defined by
$\pi_{-i}(\sigma)=\sigma_{-i}$.

A mixed strategy $\sigma_i^{*}$ of player $i$ is a best response to the
\gameenv{i} $\sigma_{-i}$ if for every mixed strategy $\sigma_i$ of player $i$,
we have $u_i(\sigma_i^{*},\sigma_{-i})\geq u_i(\sigma_i,\sigma_{-i})$.  A
strategy profile $\sigma^{*}$ is a \emph{Nash equilibrium} if for each player
$i$, the mixed strategy $\sigma_i^{*}$ is a best response to $\sigma_{-i}^{*}$.  

We can rewrite the expected payoff to player $i$ as follows:
\begin{eqnarray*}
u_i(\sigma_i,\sigma_{-i})
&=&\sum_{s\in S} u_i(s_i,s_{-i})\sigma_1(s_1)\cdots\sigma_N(s_N)\\
&=&\sum_{s_i\in S_i}\sum_{s_{-i}\in S_{-i}}u_i(s_i,s_{-i})
\sigma_1(s_1)\cdots\sigma_N(s_N)\\
&=&\sum_{s_i\in S_i}\sigma_i(s_i)
\sum_{s_{-i}\in S_{-i}}u_i(s_i,s_{-i})
\sigma_1(s_1)\cdots\sigma_{i-1}(s_{i-1})
\sigma_{i+1}(s_{i+1})\cdots\sigma_N(s_N)\\
&=&\sum_{j=0}^{d_i}\sigma_{ij}\sum_{s_{-i}\in S_i}u_i(s_{ij},\sigma_{-i}).
\end{eqnarray*}
We see that the function $u_i(s_{ij},\sigma_{-i})$ is a polynomial with
real-valued coefficients.  For each 
$s_{-i}=
(s_{1{j_1}},\ldots,s_{(i-1){j_{i-1}}},s_{(i+1){j_{i+1}}},\ldots,s_{N{j_N}})\in
S_{-i}$, the contribution to $u_i(s_{ij},\sigma_{-i})$ from the outcome
$(s_{ij},s_{-i})$ is the squarefree monomial
$$\sigma_{1{j_1}}\cdots\sigma_{(i-1){j_{i-1}}}
\sigma_{(i+1){j_{i+1}}}\cdots\sigma_{N{j_N}}$$
(which may be interpreted as the conditional probability that the outcome will
occur, given that player $i$ chooses pure strategy $s_{ij}$), multiplied by the
real-valued coefficient $u_i(s_{ij},s_{-i})$ (the payoff to player $i$ for that
outcome).  The following proposition describes the system of polynomial
equations whose nonnegative solutions give the Nash equilibria.

\begin{proposition}
\label{bigsystem}
The sequence of real numbers
$$(\sigma_{10},\ldots,\sigma_{1{d_1}},\ldots,
\sigma_{N0},\ldots,\sigma_{N{d_N}})$$
constitutes a Nash equilibrium if and only if for some sequence of real
numbers 
$$(v_{10},\ldots,v_{1{d_1}},\ldots,v_{N0},\ldots,v_{N{d_N}}),$$ 
the $\sigma_{ij}$'s and
$v_{ij}$'s satisfy the following system $(*)$ of $2(D+N)$ 
polynomial equations in $2(D+N)$ unknowns:
\begin{eqnarray*}
u_i(s_{ij},\sigma_{-i})+v_{ij}&=&u_i(s_{i0},\sigma_{-i})+v_{i0}
\quad\kern0.7em\hbox{for
each $i\in I$ and for $j=1,\ldots,d_i$},\\
\sigma_{ij}v_{ij}&=&0
\qquad\qquad\qquad\qquad\hbox{for each $i\in I$ and for $j=0,\ldots,d_i$},\\
\sum_{j=0}^{d_i}\sigma_{ij}&=&1
\qquad\qquad\qquad\qquad\hbox{for each $i\in I$},\\
\end{eqnarray*}
and all the $\sigma_{ij}$'s and $v_{ij}$'s are nonnegative.
\end{proposition}

\begin{proof}
$\Longrightarrow$: Suppose $\sigma$ is a Nash equilibrium.  Certainly all the
$\sigma_{ij}$'s are nonnegative and $\sum_{j=0}^{d_i}\sigma_{ij}=1$ for each
$i\in I$.  For each $i\in I$ and $j=0,\ldots,d_i$, let
$v_{ij}=u_i(\sigma)-u_i(s_{ij},\sigma_{-i})$.  Then $v_{ij}\geq 0$ since
$\sigma_i$ is a best response to $\sigma_{-i}$.  The quantities
$u_i(s_{ij},\sigma_{-i})+v_{ij}$ are all equal to $u_i(\sigma)$ and hence to
each other.  

It remains to show that $\sigma_{ij}v_{ij}=0$.  If $\sigma_{ij}=1$, then
$u_i(\sigma)=u_i(s_{ij},\sigma_{-i})$ and $v_{ij}=0$, so we may assume
$\sigma_{ij}<1$.  Define $\sigma_i^{*}$ by $\sigma^{*}_{ij}=0$
and 
$$\sigma^{*}_{il}
=\frac{\sigma_{il}}{1-\sigma_{ij}}\qquad\hbox{for $l\neq j$}$$
Clearly $\sigma^{*}_{il}\geq 0$, and 
$$\sum_{l=0}^{d_i}\sigma^{*}_{il}
=\frac{\sum_{l=0\atop l\neq
j}^{d_i}\sigma_{il}}{1-\sigma_{ij}}=\frac{1-\sigma_{ij}}{1-\sigma_{ij}}=1.$$
Then
\begin{eqnarray*}
u_i(\sigma_i^{*},\sigma_{-i})
&=&\sum_{l=0}^{d_i}\sigma^{*}_{il}u_i(s_{il},\sigma_{-i})\\
&=&\sum_{l=0\atop l\neq j}^{d_i}
\frac{\sigma_{il}}{1-\sigma_{ij}}u_i(s_{il},\sigma_{-i})\\
&=&\left(\sum_{l=0}^{d_i}
\frac{\sigma_{il}}{1-\sigma_{ij}}u_i(s_{il},\sigma_{-i})\right)
-\frac{\sigma_{ij}}{1-\sigma_{ij}}u_i(s_{ij},\sigma_{-i})\\
&=&\frac{1}{1-\sigma_{ij}}u_i(\sigma)
-\frac{\sigma_{ij}}{1-\sigma_{ij}}u_i(s_{ij},\sigma_{-i})\\
&=&u_i(\sigma)
+\frac{\sigma_{ij}}{1-\sigma_{ij}}
\left(u_i(\sigma)-u_i(s_{ij},\sigma_{-i})\right)\\
&=&u_i(\sigma)+\frac{\sigma_{ij}v_{ij}}{1-\sigma_{ij}},
\end{eqnarray*}
so $\sigma_{ij}v_{ij}=0$ since $\sigma_i$ is a best response to $\sigma_{-i}$.

\noindent
$\Longleftarrow$: Suppose the $\sigma_{ij}$'s and $v_{ij}$'s are nonnegative
and satisfy the polynomial system.  Let $\sigma$ be the strategy profile
defined by $\sigma_i(s_{ij})=\sigma_{ij}$.  Fix a player $i$ and suppose
$\sigma_i'$ is a mixed strategy of player $i$.  Then 
$\sum_{j=0}^{d_i}\left(\sigma_{ij}'-\sigma_{ij}\right)=0$, so there must be
some $j$ for which $\sigma_{ij}'-\sigma_{ij}\geq 0$.  Without loss of
generality suppose that $\sigma_{i0}'-\sigma_{i0}\geq 0$.  Then
\begin{eqnarray*}
u_i(\sigma)-u_i(\sigma_i',\sigma_{-i})
&=&\sum_{j=0}^{d_i}(\sigma_{ij}-\sigma_{ij}')u_i(s_{ij},\sigma_{-i})\\
&=&\sum_{j=0}^{d_i}(\sigma_{ij}-\sigma_{ij}')
\left(u_i(s_{ij},\sigma_{-i})-u_i(s_{i0},\sigma_{-i})\right)\\
&=&\sum_{j=1}^{d_i}(\sigma_{ij}-\sigma_{ij}')(v_{i0}-v_{ij})\\
&=&\left((1-\sigma_{i0})-(1-\sigma_{i0}')\right)v_{i0}
+\sum_{j=1}^{d_i}\sigma_{ij}'v_{ij}\\
&=&(\sigma_{i0}'-\sigma_{i0})v_{i0}+\sum_{j=1}^{d_i}\sigma_{ij}'v_{ij}\geq0.\\
\end{eqnarray*}
Thus $\sigma$ is a best response to $\sigma_{-i}$ for each $i$, and $\sigma$ is
a Nash equilibrium.
\end{proof}

We call the $v_{ij}$'s \emph{complementary slack variables}.  If $v_{ij}>0$,
that is, the payoff $u_i(s_{ij},\sigma_{-i})$ to player $i$ for strategy
$s_{ij}$ is strictly less than the equilibrium payoff $u_i(\sigma)$, then
$\sigma_{ij}=0$, that is strategy $s_{ij}$ cannot be a component of the mixed
strategy of player $i$.  Conversely, if $\sigma_{ij}>0$ (which must hold for
some $j$), then the payoff to player $i$ of pure strategy $s_{ij}$ must equal
the equilibrium payoff.

\section{What Kind of Geometric Object Is The Set of Nash Equilibria?}
\noindent
Thus the set of Nash equilibria is the set of real solutions to a system of
polynomial equations and inequalities with real coefficients.  Inequalities can
only be defined over ordered sets; the field of complex numbers, for example,
is not ordered.  A set of real points given by a system of polynomial equations
and inequalities is called a {\em semialgebraic variety}, and the special case
when the system does not involve inequalities is called a {\em real algebraic
variety}.  Thus the set of Nash equilibria of a game is a semialgebraic
variety.  \emph{Real algebraic geometry} is the study of real algebraic
varieties and semialgebraic varieties.  This area of algebraic geometry has
many features of special interest (quite apart from its usefulness in
applications).  (For example, over the real numbers
any system of equations $f_1=0,\ldots,f_n=0$ is equivalent to a single equation
$f_1^2+f_2^2+\cdots+f_n^2=0$.)  It so happens that Nash's contribution to real
algebraic geometry was seminal, although he did not relate it to game theory.
The two main references for real algebraic geometry are \cite{BochnakCosteRoy} 
and
\cite{BasuPollackRoy}.  

The set of points in the plane with $x\geq 0$ and $y\geq 0$ is the nonnegative
\emph{quadrant}, the set of points in 3-space with $x\geq 0$, $y\geq 0$, and
$z\geq 0$ is the nonnegative \emph{octant}, and similarly the set of points in
$\R^n$ all of whose coordinates are nonnegative is the nonnegative
\emph{orthant}.  In this case, the inequalities state simply that we are
interested in those solutions to the polynomial system which lie in the
nonnegative orthant.

In the system $(*)$ of polynomial equations, we can substitute each
$\sigma_{ij}$ with $\rho_{ij}^2$ and each $v_{ij}$ with $r_{ij}^2$, where the
$\rho_{ij}$'s and the $r_{ij}$'s are new unknowns.  This induces a new system
$(**)$ of polynomial equations in the $\rho_{ij}$'s and $r_{ij}$'s.  (For
example, $\rho_{ij}^2r_{ij}^2=0$, $\sum_{j=0}^{d_i}\rho_{ij}^2=1$, and so
forth.)  Each real-valued solution to $(**)$ corresponds to a Nash equilibrium, 
since $\sigma_{ij}=\rho_{ij}^2$ and $v_{ij}=r_{ij}^2$ automatically satisfy the
nonnegativity constraints.  However, there may be up to $2^{D+N}$ solutions to
$(**)$ for each Nash equilibrium, since if $\sigma_{ij}$ is positive we can set
each $\rho_{ij}$ equal to either its positive or negative square root, and
similarly for $v_{ij}$.  

A \emph{transformation} of a set $X$ is a 1-1 correspondence of $X$ with
itself.  A \emph{group of transformations} is a set $G$ of transformations such
that for each transformation in $G$, its inverse transformation is also in
$G$, and for each pair of transformations in $G$, their composition is
also in $G$.  We consider the composition of two transformations in $G$ to be
their product in $G$.  (Note that this product is not necessarily commutative.)
We say that the group of transformations \emph{acts} on the
set $X$ which it transforms.  Given any set $T$ of transformations
$\{g_1,\ldots,g_n\}$, we can form the group \emph{generated} by the
\emph{generators} $g_1,\ldots,g_n$ by taking all products of elements of $T$
and their inverses.  An \emph{orbit} of the group action is the set of images
of a single point in $X$ under all the transformations in $G$.  Every point in
an orbit is the image of every other point under some transformation in $G$.
To form the \emph{quotient} of $X$ by $G$, we can take one point from each
orbit.

Each $(n-1)$-dimensional hyperplane in $\R^n$ defines a special transformation,
the \emph{reflection} which takes each point in $\R^n$ to its opposite point on
the other side of the hyperplane.  In particular, the transformation of
$\R^{2(D+N)}$ which takes $\rho_{ij}$ into $-\rho_{ij}$ and leaves all other
coordinates unchanged is a \emph{coordinate reflection}, as is the one which
takes $r_{ij}$ into $-r_{ij}$.  A group of transformations generated by
reflections is called a \emph{reflection group}.  In particular, let $G$ be the
group of transformations generated by the aforementioned coordinate
reflections.  Then for any transformation $g\in G$, for any real-valued solution
$(\rho,r)$ of $(**)$, its image $g(\rho,r)$ also satisfies $(**)$.  Indeed, $g$
simply changes the signs of some of the coordinates of $(\rho,r)$.  Thus the
set of real-valued solutions to $(**)$ is a \emph{symmetric} real algebraic
variety $V$, with $G$ as its group of symmetries.  

We can take the quotient of $V$ by $G$ by, for example, considering only those
points of $V$ lying in the nonnegative orthant.  The group $G$ takes this
orthant to each other orthant.  There is exactly one point of $V$ in the
nonnegative orthant for each Nash equilibrium, obtained by taking the
nonnegative square root of each coordinate of the Nash equilibrium.  So the set
of Nash equilibria is the quotient of a real algebraic variety by a reflection
group.  In fact, any semialgebraic variety defined by a system of polynomial
equations and inequalities in which none of the inequalities are strict, is
also the quotient of a real algebraic variety by a reflection group.  For each
inequality constraint $f(x_1,\ldots,x_n)\geq0$, we can define a new variable
$v$ and change the inequality constraint into two constraints
$f(x_1,\ldots,x_n)-v=0$ and $v\geq0$, then proceed as above.

\section{Eliminating The Complementary Slack Variables}
\noindent
Fix a player $i\in I$, a pure strategy $s_{ij}\in S_i$ of
player $i$, and a pure strategy $s_{k{l_k}}\in S_k$ for every other player
$k\in I-\{i\}$, giving a pure \gameenv{i} $s_{-i}\in S_{-i}$.  Then the system
$(*)$ implies that 
$$u_i(s_{ij},\sigma_{-i})+v_{ij}
=u_i(s_{i0},\sigma_{-i})+v_{i0}.$$
Multiplying through by 
$$\prod_{k\in I-\{i\}}\prod_{m_k=0\atop m_k\neq l_k}^{d_k}v_{k{m_k}},$$
which we denote by $p(s_{-i})$,
and using that $\sigma_{k{m_k}}v_{k{m_k}}=0$, we find that
\begin{eqnarray*}
\left(u_i(s_{ij},s_{-i})\prod_{k\in I-\{i\}}\sigma_{k{l_k}}
+v_{ij}\right)
\prod_{k\in I-\{i\}}\prod_{m_k=0\atop m_k\neq l_k}^{d_k}v_{k{m_k}}\\
\qquad=
\left(u_i(s_{i0},s_{-i})\prod_{k\in I-\{i\}}\sigma_{k{l_k}}
+v_{i0}\right)
\prod_{k\in I-\{i\}}\prod_{m_k=0\atop m_k\neq l_k}^{d_k}v_{k{m_k}}.\\
\end{eqnarray*}
Every other term in $u_i(s_{ij},\sigma_{-i})$ is killed by one or more of
the $v_{k{m_k}}$'s, and similarly for $u_i(s_{i0},\sigma_{-i}$.  
Next we substitute
$1-\sum_{{m_k}=0\atop{m_k}\neq{l_k}}^{d_k}\sigma_{k{m_k}}$ for
$\sigma_{k{l_k}}$ in each factor of the product
$\prod_{k\in I-\{i\}}\sigma_{k{l_k}}$, and note that every
nonconstant term in the resulting polynomial is also killed by one or more of
the $v_{k{m_k}}$'s.  Finally, we have
\[\left(v_{ij}-v_{i0}
+u_i(s_{ij},s_{-i})-u_i(s_{i0},s_{-i})\right)
\prod_{k\in I-\{i\}}\prod_{m_k=0\atop m_k\neq l_k}^{d_k}v_{k{m_k}}=0,\]
an equation involving only the complementary slack variables.  Clearly, for
this equation to hold, either
$v_{ij}-v_{i0}+u_i(s_{ij},s_{-i})-u_i(s_{i0},s_{-i})$ or one of
the $v_{k{m_k}}$'s must vanish.  

In this way we obtain a system of equations in the $v_{ij}$'s, which we denote
by \V{}.  Denote by \V{_i} the subsystem corresponding to fixing player $i$
above, and by \V{_{i,{s_{-i}}}} the subsystem corresponding to fixing player $i$
and the pure \gameenv{i} $s_{-i}$.

If $p(s_{-i})\neq0$ for some $s_{-i}\in S_{-i}$, then the system of equations
\V{_{i,{s_{-i}}}} reduces to
$v_{ij}=v_{i0}+u_i(s_{i0},s_{-i})-u_i(s_{ij},s_{-i})$ for each
$j=1,\ldots,d_i$.  The solutions to this system along with the nonnegativity
constraints $v_{ij}\geq 0$ lie along a ray parametrized by
$v_{i0}$, with
\[v_{i0}\geq
\max_{j=0}^{d_i}\left(u_i(s_{ij},s_{-i})-u_i(s_{i0},s_{-i})\right).\]
Only for those $j=0,\ldots,d_1$ at which
$u_i(s_{ij},s_{-i})-u_i(s_{i0},s_{-i})$ attains its maximum can $v_{ij}$ ever
vanish; the rest of the $v_{ij}$'s must be positive.  But since
$\sum_{j=0}^{d_i}\sigma_{ij}=1$ with $\sigma_{ij}\geq 0$, at least one
$\sigma_{ij}$ must be positive, and thus since $\sigma_{ij}v_{ij}=0$ for all
$j$, at least one $v_{ij}$ must vanish.  Therefore we have a unique solution,
the point with
$v_{i0}=\max_{j=0}^{d_i}\left(u_i(s_{ij},s_{-i})-u_i(s_{i0},s_{-i})\right)$.
Let us denote this point by $q(s_{-i})\in\R^{d_i+1}_{\geq 0}$.

For generic games, the point $q(s_{-i})$ corresponding to
$s_{-i}$ will be different from the point $q(s'_{-i})$
corresponding to $s'_{-i}$ for any other $s'_{-i}\in S_{-i}$.  
Thus the system \V{_i} reduces to the
following $|S_{-i}|+1$ alternatives: either 
\[p(s_{-i})=0\qquad\hbox{for all $s_{-i}\in S_{-i}$},\eqno{({\rm\bf V}^p_i)}\]
or for some $s^{*}_{-i}\in S_{-i}$, we have the system \V{^{s^{*}_{-i}}_i}:
\[p(s_{-i})=0\qquad\hbox{for all $s_{-i}\in S_{-i}-\{s^{*}_{-i}\}$},
\eqno{({\rm\bf V}^{s^{*}_{-i}p}_i)}\]
and
\[(v_{i0},\ldots,v_{i{d_i}})=q(s^{*}_{-i}).
\eqno{({\rm\bf V}^{s^{*}_{-i}q}_i)}\]
For generic games, exactly one of the $v_{ij}$'s will vanish at the point
$q(s^{*}_{-i})$.

We now introduce some notation from logic.  The symbol $\wedge$ means ``and''
and the symbol $\vee$ means ``or'' (in the Boolean sense).  Suppose $\varphi_0$
and $\varphi_1$ are Boolean expressions. Then the expression
$\varphi_0\wedge\varphi_1$, meaning ``$\varphi_0$ and $\varphi_1$'', is the
\emph{conjunction} of $\varphi_0$ and $\varphi_1$.  The expression
$\varphi_0\vee\varphi_1$, meaning ``$\varphi_0$ or $\varphi_1$'', is the
\emph{disjunction} of $\varphi_0$ and $\varphi_1$.  Suppose
$\varphi_0,\ldots,\varphi_n$ are Boolean expressions.  Then
\(\bigwedge_{i=0}^n\varphi_i\)
denotes the conjunction of $\varphi_0,\ldots,\varphi_n$, and
\(\bigvee_{i=0}^n\varphi_i\)
denotes the disjunction of $\varphi_0,\ldots,\varphi_n$.  The distributive law
holds for conjunction and disjunction just as it holds for multiplication and
addition.

\begin{lemma}
\label{allprodzero}
The solutions of the system \V{^p_i} are given by
\[\bigvee_{k=1\atop k\neq i}^N \bigvee_{j=0}^{d_k} 
\bigvee_{l=1\atop l>j}^{d_k}
\left(v_{kj}=0\right)\wedge\left(v_{kl}=0\right).\]
\end{lemma}

\begin{proof}
We proceed by induction on the number $N$ of players.  First suppose $N=2$.
Without loss of generality, fix $i=2$.  Then each $s_{-i}$ is $s_{1j}$ for some
$j$.  The equation \V{^p_{2,s_{1j}}} is
\(p(s_{1j})=\prod_{l=0\atop l\neq j}^{d_1}v_{1l}=0\),
which holds if and only if 
\(\bigvee_{l=0\atop l\neq j}^{d_1}\left(v_{1l}=0\right)\).
So the system \V{^p_2} holds if and only if
\(\bigwedge_{j=0}^{d_1}\bigvee_{l=0\atop l\neq j}^{d_1}\left(v_{1l}=0\right)\).
We show by induction on $d_1$ that this expression is equal to
\(\bigvee_{j=0}^{d_1}\bigvee_{l=1\atop l>j}^{d_1}
\left(v_{1j}=0\right)\wedge\left(v_{1l}=0\right)\)
which equals
\(\bigvee_{k=1\atop k\neq 2}^2
\bigvee_{j=0}^{d_k}\bigvee_{l=1\atop l>j}^{d_k}
\left(v_{kj}=0\right)\wedge\left(v_{kl}=0\right)\).
For $d_1=1$, 
\begin{eqnarray*}
\bigwedge_{j=0}^{1}\bigvee_{l=0\atop l\neq j}^{1}\left(v_{1l}=0\right)
&=&\bigwedge_{j=0}^1\left(v_{1(1-j)}=0\right)\\
&=&\left(v_{11}=0\right)\wedge\left(v_{10}=0\right)\\
&=&\bigvee_{j=0}^{1}\bigvee_{l=1\atop l>j}
\left(v_{1j}=0\right)\wedge\left(v_{1l}=0\right)
\end{eqnarray*}
Now assume that 
\(\bigwedge_{j=1}^{d_1}\bigvee_{l=1\atop l\neq j}^{d_1}\left(v_{1l}=0\right)
=\bigvee_{j=1}^{d_1}\bigvee_{l=2\atop l>j}^{d_1}
\left(v_{1j}=0\right)\wedge\left(v_{1l}=0\right)\).
Then
\begin{eqnarray*}
\bigwedge_{j=0}^{d_1}\bigvee_{l=0\atop l\neq j}^{d_1}\left(v_{1l}=0\right)
&=&\biggl(\bigvee_{l=1}^{d_1}\left(v_{1l}=0\right)\biggr)\wedge
\bigwedge_{j=1}^{d_1}\bigvee_{l=0\atop l\neq j}^{d_1}\left(v_{1l}=0\right)\\
&=&\biggl(\bigvee_{l=1}^{d_1}\left(v_{1l}=0\right)\biggr)\wedge
\bigwedge_{j=1}^{d_1}\Bigl(\left(v_{10}=0\right)\vee
\bigvee_{l=1\atop l\neq j}^{d_1}\left(v_{1l}=0\right)\Bigr)\\
&=&\biggl(\bigvee_{l=1}^{d_1}\left(v_{1l}=0\right)\biggr)\wedge
\Bigl(\left(v_{10}=0\right)\vee
\bigwedge_{j=1}^{d_1}\bigvee_{l=1\atop l\neq j}^{d_1}
\left(v_{1l}=0\right)\Bigr)\\
&=&\biggl(\bigvee_{l=1}^{d_1}\left(v_{1l}=0\right)\biggr)\wedge
\Bigl(\left(v_{10}=0\right)\vee
\bigvee_{j=1}^{d_1}\bigvee_{l=2\atop l>j}^{d_1}
\left(v_{1j}=0\right)\wedge\left(v_{1l}=0\right)\Bigr)\\
&=&\biggl(\bigvee_{l=1}^{d_1}
\left(v_{10}=0\right)\wedge\left(v_{1l}=0\right)\biggr)\\
&&\qquad\vee\Biggl(
\biggl(\bigvee_{l=1}^{d_1}\left(v_{1l}=0\right)\biggr)\wedge
\bigvee_{j=1}^{d_1}
\biggl(
\left(v_{1j}=0\right)\wedge
\Bigl(\bigvee_{l=1\atop l>j}^{d_1}\left(v_{1l}=0\right)\Bigl)\biggr)\Biggr)\\
&=&\biggl(\bigvee_{j=0}^0\bigvee_{l=1\atop l>j}^{d_1}
\left(v_{1j}=0\right)\wedge\left(v_{1l}=0\right)\biggr)\\
&&\qquad\vee\Biggl(
\biggl(\bigvee_{j=1}^{d_1}\left(v_{1j}=0\right)\biggr)\wedge
\bigvee_{j=1}^{d_1}
\biggl(
\left(v_{1j}=0\right)\wedge
\Bigl(\bigvee_{l=1\atop l>j}^{d_1}\left(v_{1l}=0\right)\Bigl)\biggr)\Biggr)\\
&=&\biggl(\bigvee_{j=0}^0\bigvee_{l=1\atop l>j}^{d_1}
\left(v_{1j}=0\right)\wedge\left(v_{1l}=0\right)\biggr)\vee
\bigvee_{j=1}^{d_1}
\biggl(
\left(v_{1j}=0\right)\wedge
\Bigl(\bigvee_{l=1\atop l>j}^{d_1}\left(v_{1l}=0\right)\Bigl)\biggr)\\
&=&\bigvee_{j=0}^{d_i}\bigvee_{l=1\atop l>j}^{d_1}
\left(v_{1j}=0\right)\wedge\left(v_{1l}=0\right).\\
\end{eqnarray*}

\noindent
Now assume the lemma holds for any number of players less than $N$.  Fix $i\in
I-\{N\}$.  Suppose that 
\(\bigvee_{k=1\atop k\neq i}^N \bigvee_{j=0}^{d_k} 
\bigvee_{l=1\atop l>j}^{d_k}
\left(v_{kj}=0\right)\wedge\left(v_{kl}=0\right)\)
holds.  Fix $n\in I-\{i\}$ and $j,m\in S_n$ with $j\neq m$ such that $v_{nj}=0$
and $v_{nm}=0$.  
Consider a pure strategy $s_{k{l_k}}\in S_k$ for
every player $k\in I-\{1\}$, giving a pure \gameenv{1} 
$s_{-i}\in S_{-i}$.  
The equation \V{^p_{i,s_{-i}}} is
\[p(s_{-i})=\prod_{k\in I-\{i\}}\prod_{m_k=0\atop m_k\neq l_k}^{d_k}v_{k{m_k}}
=0.\]
Either $l_n=j$, in which case the product vanishes since it includes $v_{nm}$;
or $l_n=m$, in which case the product vanishes since it includes $v_{nj}$; or
$s_{n{l_n}}\in S_n-\{s_{nj},s_{nm}\}$, in which case the product vanishes since
it includes both $v_{nj}$ and $v_{nm}$.  Thus the equation \V{^p_{i,s_{-i}}}
holds.  So the system \V{^p_i} is satisfied.

Conversely, suppose the system \V{^p_i} holds.  If there exist $j,m\in S_N$
with $j\neq m$ such that $v_{Nj}=0$ and $v_{Nm}=0$, then
\(\bigvee_{k=1\atop k\neq i}^N \bigvee_{j=0}^{d_k} 
\bigvee_{l=1\atop l>j}^{d_k}
\left(v_{kj}=0\right)\wedge\left(v_{kl}=0\right)\)
holds, so assume there do not exist such $j,m\in S_N$.  Then there exists 
$j\in S_N$ such that $v_{Nm}\neq0$ for all $s_{Nm}\in S_N-\{s_{Nj}\}$.
Write $S'_{-i}=\prod_{k\in I-\{i,N\}}S_k$.  Consider all elements of $S_{-i}$
given by pure strategies $s_{k{l_k}}\in S_k$ for $k\in I-\{i,N\}$ and
$s_{N{l_N}}=s_{Nj}\in S_N$.  Thus $s_{-i}$ has the form
$s_{-i}=(s'_{-i},s_{Nj})$ for some $s'_{i}\in S'_{i}$.
The equation \V{^p_{i,s_{-i}}} is
\[p(s_{-i})
=\prod_{k\in I-\{i\}}\prod_{m_k=0\atop m_k\neq l_k}^{d_k}v_{k{m_k}}
=\Biggl(\prod_{k\in I-\{i,N\}}\prod_{m_k=0\atop m_k\neq l_k}^{d_k}
v_{k{m_k}}\Biggr)
\prod_{m=0\atop m\neq j}^{d_N}v_{Nm}
=0.\]
Since $v_{Nm}\neq0$ for all $s_{Nm}\in S_N-\{s_{Nj}\}$, this equation is
equivalent to the equation
\(p(s'_{-i})=\prod_{k\in I-\{i,N\}}\prod_{m_k=0\atop m_k\neq
l_k}^{d_k}v_{k{m_k}}=0\).
By the induction hypothesis, these equations imply that
\(\bigvee_{k=1\atop k\neq i}^{N-1} \bigvee_{j=0}^{d_k} 
\bigvee_{l=1\atop l>j}^{d_k}
\left(v_{kj}=0\right)\wedge\left(v_{kl}=0\right)\)
holds.  Thus
\(\bigvee_{k=1\atop k\neq i}^N \bigvee_{j=0}^{d_k} 
\bigvee_{l=1\atop l>j}^{d_k}
\left(v_{kj}=0\right)\wedge\left(v_{kl}=0\right)\)
holds.
\end{proof}

\begin{lemma}
\label{allbutoneprodzero}
Fix $i\in I$ and fix pure strategies $s_{k{l_k}}\in S_k$ for each $k\in
I-\{i\}$, defining a \gameenv{i} $s^{*}_{-i}\in S_{-i}$.  Then 
the solutions of the system \V{^{s^{*}_{-i}p}_i} are given by
\[\bigwedge_{k=1\atop k\neq i}^N \left(v_{k{l_k}}=0\right)\vee
\bigvee_{k=1\atop k\neq i}^N \bigvee_{j=0}^{d_k} 
\bigvee_{l=1\atop l>j}^{d_k}
\left(v_{kj}=0\right)\wedge\left(v_{kl}=0\right).\]
\end{lemma}

\begin{proof}
Suppose
\(\bigwedge_{k=1\atop k\neq i}^N \left(v_{k{l_k}}=0\right)\vee
\bigvee_{k=1\atop k\neq i}^N \bigvee_{j=0}^{d_k} 
\bigvee_{l=1\atop l>j}^{d_k}
\left(v_{kj}=0\right)\wedge\left(v_{kl}=0\right)\)
holds.  If
\[\bigvee_{k=1\atop k\neq i}^N \bigvee_{j=0}^{d_k} 
\bigvee_{l=1\atop l>j}^{d_k}
\left(v_{kj}=0\right)\wedge\left(v_{kl}=0\right)\]
holds, then by Lemma \ref{allprodzero}, the system \V{^p_i} holds so \emph{a
fortiori} the system \V{^{s^{*}_{-i}p}_i} holds.  Suppose 
\(\bigwedge_{k=1\atop k\neq i}^N \left(v_{k{l_k}}=0\right)\) holds.
Let $s_{-i}\in S_{-i}-\{s^{*}_{-i}\}$ with components $s_{k{n_k}}\in S_k$.  Then
\[p(s_{-i})
=\prod_{k=1\atop k\neq i}^N\prod_{m_k=0\atop m_k\neq n_k}^{d_k}v_{k{m_k}}.\]
Since $s_{-i}\neq s^{*}_{-i}$, there is some $k\in I-\{i\}$ such that 
$n_k\neq l_k$.  Then $v_{k{l_k}}$ appears in the product, and so
$p(s_{-i})=0$.  So the system \V{^{s^{*}_{-i}p}_i} holds.

Conversely, suppose 
\(\bigwedge_{k=1\atop k\neq i}^N \left(v_{k{l_k}}=0\right)\vee
\bigvee_{k=1\atop k\neq i}^N \bigvee_{j=0}^{d_k} 
\bigvee_{l=1\atop l>j}^{d_k}
\left(v_{kj}=0\right)\wedge\left(v_{kl}=0\right)\)
does not hold.  
Since
\(\bigwedge_{k=1\atop k\neq i}^N \left(v_{k{l_k}}=0\right)\)
does not hold, there is some $k'\in I-\{i\}$ such that $v_{k'{l_{k'}}}\neq0$.
Also since
\(\bigvee_{k=1\atop k\neq i}^N \bigvee_{j=0}^{d_k} 
\bigvee_{l=1\atop l>j}^{d_k}
\left(v_{kj}=0\right)\wedge\left(v_{kl}=0\right)\)
does not hold, for each $k\in I-\{i\}$, 
either $v_{kj}\neq0$ for all $s_{kj}\in S_k$, or $v_{k{n_k}}=0$ for a single
$s_{k{n_k}}\in S_k$ and 
$v_{k{m_k}}\neq0$ for all $m_k$ with $s_{k{m_k}}\in S_k-\{s_{k{n_k}}\}$.  
In either case, there is some $n_k$ with $s_{k{n_k}}\in S_k$ such that
$v_{k{m_k}}\neq0$ for all $m_k$ with $s_{k{m_k}}\in S_k-\{s_{k{n_k}}\}$.  
In particular, for $k'$ we can choose $n_{k'}\neq l_{k'}$, because in the
former case, any $n_{k'}$ other than $l_{k'}$ will do, and in the latter case,
the single $n_{k'}$ with $v_{k{n_{k'}}}=0$ must not be equal to $l_{k'}$, since
$v_{k'{l_{k'}}}\neq0$.
Define $s_{-i}\in S_{-i}$ with components $s_{k{n_k}}\in S_k$
such that $v_{k{m_k}}\neq0$ for all $m_k$ with 
$s_{k{m_k}}\in S_k-\{s_{k{n_k}}\}$,
for all $k\in I-\{i,k'\}$, and $n_{k'}\neq l_{k'}$.  Then
\[p(s_{-i})
=\prod_{k=1\atop k\neq i}^N\prod_{m_k=0\atop m_k\neq n_k}^{d_k}v_{k{m_k}}\neq0\]
by the choice of the $n_k$'s.  Since $n_{k'}\neq l_{k'}$, we have 
$s_{-i}\in S_{-i}-\{s^{*}_{-i}\}$.  
So \V{^{s^{*}_{-i}p}_{i,s_{-i}}} does not hold and the system
\V{^{s^{*}_{-i}p}_i} does not hold.
\end{proof}

\begin{corollary}
\label{purestrictorallprodzero}
The system \V{_i} holds if and only if either
\[\bigvee_{k=1\atop k\neq i}^N \bigvee_{j=0}^{d_k} 
\bigvee_{l=1\atop l>j}^{d_k}
\left(v_{kj}=0\right)\wedge\left(v_{kl}=0\right),\]
or for some pure strategies $s_{kl}\in S_{k}$ for each $k\in I-\{i\}$ defining
a \gameenv{i} $s^{*}_{-i}\in S_{-i}$ , we have
\[(v_{i0},\ldots,v_{i{d_i}})=q(s^{*}_{-i}),\]
\[v_{kl_k}=0\qquad\hbox{for each $k\in I-\{i\}$},\]
and
\[v_{km_k}\neq0\qquad\hbox{for all $m_k\neq l_k$, for each $k\in I-\{i\}$}.\]
\end{corollary}

\begin{proof}
This follows from our characterization of \V{_i}, Lemma \ref{allprodzero} and
Lemma \ref{allbutoneprodzero}.  We impose the condition
\[v_{k{m_k}}\neq0\qquad\hbox{for all $m_k\neq l_k$, for each $k\in I-\{i\}$}\]
in the latter alternative because otherwise, the former alternative holds.
\end{proof}

Suppose we have a solution of the system \V{} in which the latter alternative
holds for some $i\in I$.  As noted earlier, for a generic game, exactly one
$v_{ij}=0$ (the one for which $u_i(s_{ij},s_{-i})-u_i(s_{i0},s_{-i})$ attains
its maximum, i.e., the one for which $u_i(s_{ij},s_{-i})$ is maximum) and the
rest are zero.  So the former alternative cannot hold for \emph{any} $k\in I$.
Define $l_i$ by $v_{i{l_i}}=0$, and let $s^{*}=(s_{i{l_i}},s^{*}_{-i})\in S$.  
For each $k\in I$, the payoff $u_k(s_{k{l_k}},s^{*}_{-k})$ must be the maximum
among $u_k(s_{kl},s^{*}_{-k})$.  Since $v_{kl}\neq 0$ for $l\neq l_k$, we must
have $\sigma_{kl}=0$ for $l\neq l_k$, so $\sigma_{k{l_k}}=1$.  That is,
$\sigma$ is the pure strategy profile $s^{*}$.  So this is the case of a pure
\emph{strict} Nash equilibrium $s^{*}$, that is, one for which the pure
strategy $s^{*}_k$ is a strictly better response to $s^{*}_{-k}$ than any other
pure strategy of $k$, for each $k\in I$.  

Note well that this is a weaker condition than that $s^{*}$ be an equilibrium
in dominant strategies.  For example, consider a game of two players who each
can take one of two actions, in which the payoff to each player is the same if
they take the same action and strictly less if they take opposite actions.
This game has two pure strict Nash equilibria, corresponding to both players
taking the same one of the two actions.  But neither of the actions is a
dominated strategy for either of the players.

Suppose the latter alternative holds for some $i\in I$ and both $v_{ij}=0$ and
$v_{ij'}=0$ for some $j\neq j'$.  Then this makes the former alternative true
for every $k\in I-\{i\}$.

Finding pure strict Nash equilibria is a combinatorial procedure which does not
require any polynomial algebra.  Therefore we do not discuss it further in this
article, but assume that we have already found all pure strict Nash equilibria
(if any exist), and are now interested in finding the other ones.

\begin{proposition}
\label{twobytwo}
Suppose $\sigma$ is a Nash equilibrium of a generic game, and $\sigma$ is not a
pure strict Nash equilibrium.  Then there are two players $i,k\in I$ with
$i\neq k$ and two pure strategies each, $s_{i{j_0}},s_{i{j_1}}\in S_i$ with
$j_0\neq j_1$ and $s_{k{l_0}},s_{k{l_k}}\in S_k$ with $l_0\neq l_1$, such that
$u_i(s_{i{j0}},\sigma_{-i})=u_i(s_{i{j1}},\sigma_{-i})=u_i(\sigma)$ and
$u_k(s_{k{l0}},\sigma_{-k})=u_k(s_{k{l1}},\sigma_{-k})=u_k(\sigma)$.
\end{proposition}

\begin{proof}
The first alternative in Corollary \ref{purestrictorallprodzero} must hold for
each $i\in I$.  Pick a player $n\in I$.  The condition \[\bigvee_{k=1\atop
k\neq n}^N \bigvee_{j=0}^{d_k} \bigvee_{l=1\atop l>j}^{d_k}
\left(v_{kj}=0\right)\wedge\left(v_{kl}=0\right)\] means there is $k\neq n$ and
pure strategies $s_{k{l_0}},s_{k{l_1}}\in S_k$ with $v_{k{l_0}}=v_{k{l_1}}=0$,
i.e., $u_k(s_{k{l_0}},\sigma_{-k}=u_k(s_{k{l_1}},\sigma_{-k})=u_k(\sigma)$.  This
makes \V{_m} hold for every player $m\in I-\{k\}$.  For $k$ itself, the
condition \(\bigvee_{i=1\atop i\neq k}^N \bigvee_{j=0}^{d_i} \bigvee_{l=1\atop
l>j}^{d_i} \left(v_{ij}=0\right)\wedge\left(v_{il}=0\right)\) means there is
$i\neq k$ and pure strategies $s_{i{j_0}},s_{i{j_1}}\in S_i$ with
$u_i(s_{i{j0}},\sigma_{-i})=u_i(s_{i{j1}},\sigma_{-i})=u_i(\sigma)$.
\end{proof}

In this case we cannot isolate $\sigma_{i{j_0}}$ and $\sigma_{i{j_1}}$ simply
by looking at the best responses for $i$, and similarly with $k$.  We must
solve the polynomial system $(*)$.

As noted before, at least one $v_{nm}$ must vanish for every $n\in I$.  In the
conditions of the proposition, the least complex case is that exactly one
$v_{nm}$ vanishes for each $n\in I-\{i,k\}$, say $v_{n{m_n}}$; that
$v_{ij}>0$ for $j\neq j_0$ and $j\neq j_1$; 
and that $v_{kl}>0$ for $k\neq k_0$ and $k\neq k_1$.
Then each player $n\in I-\{i,k\}$ executes pure strategy $s_{n{m_k}}$.
Furthermore $i$ does not execute pure strategy $s_{ij}$ with any probability
for $j\neq j_0$ and $j\neq j_1$, and $k$ does not execute pure strategy
$s_{kl}$ with any probability for $l\neq l_0$ and $l\neq l_1$.  Then the system
$(*)$ reduces to the system for a game with two players, which we renumber as 1
and 2, with two pure strategies each, which we renumber as $s_{10}$, $s_{11}$,
$s_{20}$, and $s_{21}$.  Write $u^{i}_{jl}=u_i(s_{0j},s_{1l})$.  In this case
the system $(*)$ is:
\begin{eqnarray*}
u^1_{00}\sigma_{20}+u^1_{01}\sigma_{21}
&=&u^1_{10}\sigma_{20}+u^1_{11}\sigma_{21},\\
u^2_{00}\sigma_{10}+u^2_{10}\sigma_{11}
&=&u^2_{01}\sigma_{10}+u^2_{11}\sigma_{11},\\
\sigma_{10}+\sigma_{11}&=&1,\\
\sigma_{20}+\sigma_{21}&=&1.
\end{eqnarray*}
Substituting $1-\sigma_{11}$ for $\sigma_{10}$ and $1-\sigma_{21}$ for
$\sigma_{20}$, we obtain
\begin{eqnarray*}
(u^1_{11}-u^1_{10}-u^1_{01}+u^1_{00})\sigma_{21}&=&u^1_{00}-u^1_{10},\\
(u^2_{11}-u^2_{10}-u^2_{01}+u^2_{00})\sigma_{11}&=&u^2_{00}-u^2_{01}.
\end{eqnarray*}
Notice that the equilibrium found by this system need not be totally mixed; for
instance, $u^1_{00}-u^1_{10}$ could equal zero, in which case
$\sigma_2=s_{20}$, or $u^1_{11}-u^1_{01}$ could equal zero. in which case
$\sigma_2=s_{21}$.  (If they are both zero, then player $1$ has no control over
player $1$'s own payoff and hence every mixed strategy of player $1$ is a best
response.  So the requirement that $1$ play a best response, which leads to the
first equation, does not impose any constraint on the strategy $\sigma_{21}$ of
player $2$.) Similarly, $u^2_{00}-u^2_{01}$ could equal zero, in which case
$\sigma_1=s_{10}$, or $u^2_{11}-u^2_{10}$ could equal zero, in which case
$\sigma_1=s_{11}$.  Both of these cases may even occur, in which case the root
of the system is a pure Nash equilibrium.

Any solution to the system $(*)$ induces a partition of 
$P{\buildrel {\rm def}\over=}\bigcup_{i\in I}S_i$
into two subsets, the subset $P_0$ such that $v_{ij}=0$ for all $s_{ij}\in P_0$
and the subset $P_+$ such that $v_{ij}>0$ for all $s_{ij}\in P_+$.  If we make
a choice of such a partition, then $\sigma_{ij}=0$ for all $s_{ij}\in P_+$, so
eliminating these strategies and considering the reduced game, the system $(*)$
reduces to the system \E{}:
\[u_i(s_{ij},\sigma_{-i})=u_i(s_{i0},\sigma_{-i})
\quad\kern0.7em\hbox{for
each $i\in I$ and for $j=1,\ldots,d_i$},\eqno{({\rm\bf E}_{ij})}\]
\[\sum_{j=0}^{d_i}\sigma_{ij}=1
\qquad\qquad\qquad\qquad\hbox{for each $i\in I$},\eqno{({\rm\bf E}_{i0})}\]
and all the $\sigma_{ij}$'s are nonnegative.  (A root of the polynomial
equations of the system \E{} which
does not satisfy the nonnegativity constraints is called a
\emph{quasi-equilibrium}.)  After solving the system
\E{} to find a candidate $\sigma$, we have to check that for each
strategy $s_{ij}$ of the original game which we had eliminated,
$v_{ij}=u_i(\sigma)-u_i(s_{ij},\sigma_{-i})$ is indeed nonnegative.  In that
case $\sigma$ is a solution to the original system $(*)$ and hence a Nash
equilibrium of the original game.  To find all the Nash equilibria, we can
perform this procedure for all partitions $(P_0,P_+)$ for which at least one
$s_{nm}\in P_0$ for each player $n$.  As noted above, for a generic game there
will not be a Nash equilibrium for which there is exactly one player for which 
at least two $s_{ij}$'s are in $P_0$.

From now on we will restrict our attention to solving systems of the form
\E{}.  The paper \cite{mckelveymclennan} describes this system, which
along with the constraints $\sigma_{ij}>0$ for all the $\sigma_{ij}$'s, gives
the totally mixed Nash equilibria.  The \Gambit\ software package \cite{gambit}
finds all Nash equilibria recursively,
by finding totally mixed Nash equilibria of each
reduced game by solving the corresponding system.  This algorithm for
finding all Nash equilibria is described for example in \cite{heringspeeters}
and in \cite{dattathesis}.

\section{Solving An Instance of the Polynomial System}
\noindent
For any $s\in S$, with $s_i=s_{i{j_i}}$ for each $i\in I$, we write 
\begin{eqnarray*}
u^i_{{j_1}\ldots{j_N}}
&=&u_i(s_{1{j_1}},s_{2{j_2}},\ldots,s_{(i-1){j_{i-1}}},s_{ij},
s_{(i+1){j_{i+1}}},\ldots,s_{N{j_N}})\\
&&\qquad- u_i(s_{1{j_1}},s_{2{j_2}},\ldots,s_{(i-1){j_{i-1}}},s_{i0},
s_{(i+1){j_{i+1}}},\ldots,s_{N{j_N}}).
\end{eqnarray*}
In particular $u^i_{{j_1}\ldots{j_N}}=0$ if $j_i=0$.
Then the equation \E{_{ij}} is
\[\sum_{j_1=0}^{d_1}\cdots\sum_{j_{i-1}=0}^{d_{i-1}}
\sum_{j_{i+1}=0}^{d_{i+1}}\cdots\sum_{j_N=0}^{d_N}
u^i_{{j_1}\ldots{j_{i-1}}jj_{i+1}\ldots j_N}
\sigma_{1{j_1}}\ldots\sigma_{(i-1){j_{i-1}}}
\sigma_{(i+1){j_{i+1}}}\ldots\sigma_{N{j_N}}=0.\]
How do we solve an equation like this?  We note in passing that if all the
coefficients $u^i_{{j_1}\ldots{j_{i-1}}j{j_{i+1}}\ldots {j_N}}$ had the same
sign, we would know that either all the monomials vanished (i.e.,
$\sigma_{k{l_0}}=\sigma_{k{l_1}}=0$ for some $k\in I-\{i\}$ and $1\leq
l_0<l_1\leq d_k$), or $\sigma_{ij}=0$ and $v_{ij}>v_{i0}$.  However, this
condition depends on the choice of which strategy in $S_i$ to label as
$s_{i0}$, so to check whether it ever arises would require checking the
difference between every pair 
$u_i(s_{i{j_1}},s_{-i})-u_i(s_{i{j_0}},s_{-i})$.

Well, if the equation factored, it would be easy
to solve.  That is, if we could find numbers $\mu^{(ij)}_{k{j_k}}$ for
$k\in I-\{i\}$ and $s_{k{j_k}}\in S_k$ such that
\begin{eqnarray*}
\sum_{j_1=0}^{d_1}\cdots\sum_{j_{i-1}=0}^{d_{i-1}}
\sum_{j_{i+1}=0}^{d_{i+1}}\cdots\sum_{j_N=0}^{d_N}
u^i_{{j_1}\ldots{j_{i-1}}jj_{i+1}\ldots j_N}
\sigma_{1{j_1}}\ldots\sigma_{(i-1){j_{i-1}}}
\sigma_{(i+1){j_{i+1}}}\ldots\sigma_{N{j_N}}\\
\qquad=\prod_{k\in I-\{i\}}\biggl(
\sum_{{j_k}=0}^{d_k}\mu^{(ij)}_{k{j_k}}\sigma_{k{j_k}}\biggr),
\end{eqnarray*}
then we could solve the equation by setting
$\sum_{{j_k}=0}^{d_k}\mu^{(ij)}_{k{j_k}}\sigma_{k{j_k}}$ equal to zero 
for some $k\in I-\{i\}$.

The factorization holds if and only if the coefficients of each monomial 
$$\sigma_{1{j_1}}\ldots\sigma_{(i-1){j_{i-1}}}
\sigma_{(i+1){j_{i+1}}}\ldots\sigma_{N{j_N}}$$ on both sides of the equation
are equal.  The coefficient on the left-hand side is
$u^i_{{j_1}\ldots{j_{i-1}}jj_{i+1}\ldots j_N},$
and the coefficient on the right-hand side is 
$\prod_{k\in I-\{i\}}\mu^{(ij)}_{k{j_k}}$.  Thus we have a system of
equations
$\prod_{k\in I-\{i\}}\mu^{(ij)}_{k{j_k}}=
u^i_{{j_1}\ldots{j_{i-1}}jj_{i+1}\ldots j_N}.$
Equivalently we have the system of linear equations
$\sum_{k\in I-\{i\}}\log|\mu^{(ij)}_{k{j_k}}|=
\log|u^i_{{j_1}\ldots{j_{i-1}}jj_{i+1}\ldots j_N}|$
together with sign conditions
$\prod_{k\in I-\{i\}}\sign(\mu^{(ij)}_{k{j_k}})=
\sign(u^i_{{j_1}\ldots{j_{i-1}}jj_{i+1}\ldots j_N})$.
Unfortunately this linear system is overdetermined and hence usually
inconsistent.  We have $\prod_{k=1\atop k\neq i}^N(d_i+1)$ equations in
only $\sum_{k=1\atop k\neq i}^N(d_i+1)$ unknowns $\mu^{(ij)}_{k{j_k}}$.  
So in general, equation \E{_{ij}} does not factorize.

Nevertheless, suppose all the equations \E{_{ij}} \emph{did} factorize.  How
would we solve the whole system \E{} then?  For one thing we would substitute 
$\sigma_{k0}=1-\sum_{{j_k}=1}^{d_k}\sigma_{k{j_k}}$ into each linear factor 
$\sum_{{j_k}=0}^{d_k}\mu^{(ij)}_{k{j_k}}\sigma_{k{j_k}}$ to get an affine linear
factor 
\[\mu^{(ij)}_{k0}
+\sum_{j=1}^{d_k}(\mu^{(ij)}_{k{j_k}}-\mu^{(ij)}_{k0})\sigma_{k{j_k}}.\]
(``Affine'' just means that it includes a constant term.)  We set 
$\lambda^{(ij)}_{k{j_k}}=\mu^{(ij)}_{k{j_k}}-\mu^{(ij)}_{k0}$
for $j_k=1,\ldots,d_k$ and $\lambda^{(ij)}_{k0}=\mu^{(ij)}_{k0}$.

We will now construct a particular system
in which all the equations factorize, and solve that.
For this purpose, it will be convenient to have available a \emph{totally
nonsingular matrix}.  An $m\times n$ matrix $M=(m_{ij})$ is 
\emph{totally nonsingular}
if for every $k\leq\min(m,n)$, for every subset $R\subseteq\{1,\ldots,m\}$ with
$|R|=k$
and every subset $C\subseteq\{1,\ldots,n\}$ with $|C|=k$, the $k\times k$
submatrix of $M$ given by $(m_{ij})_{i\in R\atop j\in C}$ is nonsingular.

Let $f:\N\to\R_{>0}$ be any injection of $\N$ into $\R_{>0}$.
We can use the algorithm in Figure \ref{totallynonsingularalgorithm} 
to construct a totally singular $n\times n$ matrix for any $n$.
\begin{figure}
\begin{algorithmic}
\FOR{$i=1$ to $n$}
\STATE {\it Invariant:} Every submatrix of the partially filled-in matrix is
nonsingular.
\STATE {\it Invariant:} The matrix is symmetric.
\FOR{$j=1$ to $i$}
\STATE $k\leftarrow i+j-2$
\STATE $m_{ij}\leftarrow f(k)$
\WHILE{some submatrix which includes $m_{ij}$ is singular}
\STATE $m_{ij}\leftarrow -m_{ij}$
\IF{$m_{ij}>0$}
\STATE We've already tried this value of $m_{ij}$, so try another
\STATE $k\leftarrow k+1$
\STATE $m_{ij}\leftarrow f(k)$
\ENDIF
\ENDWHILE
\STATE $m_{ji}\leftarrow m_{ij}$
\ENDFOR
\ENDFOR
\end{algorithmic}
\caption{One Possible Algorithm To Compute A Totally Nonsingular Matrix}
\label{totallynonsingularalgorithm}
\end{figure}
\noindent
We start with the $1\times 1$ matrix $(m_{11})=(f(1))$, which is clearly
totally nonsingular.  The problem of filling in a matrix while maintaining some
condition is called a \emph{matrix completion problem}.
\cite{TorregrosaJordan} shows that we can construct a totally nonsingular
matrix by filling in the entries one at a time, as in the above algorithm.  In
fact, at each stage there are only finitely many possible values of the next
entry which would violate the condition, so all we have to do is avoid those
values.  Therefore each while loop in the algorithm will always terminate.
Since the partially filled-in matrix is symmetric, and $m_{ij}$ does not
violate the condition, setting $m_{ji}=m_{ij}$ cannot violate it either, and we
can keep the matrix symmetric.  For example, here is the totally nonsingular
$6\times 6$ matrix given by the above algorithm, with $f(k)=2^{k-1}$: 

\[\left[\begin{array}{rrrrrr}
 1&   2&   4&   8&  16&   32\\
 2&  -4&  16& -32& 128& -256\\
 4&  16& -16&-128&1024& -256\\
 8& -32&-128& -64&4096& 4096\\
16& 128&1024&4096&-256& 1024\\
32&-256&-256&4096&1024&-1024\\
\end{array}\right]\]

As a matter of fact, a random matrix will be totally nonsingular with
probability one.  However, if we do use a random matrix we should check that it
is indeed totally nonsingular.  Since checking this may take a long time, it
may be useful to build a large totally nonsingular matrix once and for all and
keep it around.

Now assume we have a totally nonsingular $D\times D$ matrix $M$ with entries 
$m_{ij}$.  Define $n(i,j)=j+\sum_{k=1}^{i-1}d_k$.  So if we write the equations
\E{_{ij}} in sequence 
$(\hbox{\rm\bf E}_{11}),\ldots,(\hbox{\rm\bf E}_{1{d_1}}),
\ldots,(\hbox{\rm\bf E}_{N1}),\ldots,(\hbox{\rm\bf E}_{N{d_N}})$, then the 
$n(i,j)$th equation in the sequence is \E{_{i,j}}.  
Set $\lambda^{(ij)}_{k{l_k}}=m_{n(i,j)l_k}$ for $l_k>0$ and 
$\lambda^{(ij)}_{k0}=-1$.  This defines a particular system \startsystem{}
of equations which factorizes.

Notice that we don't use all $D^2$ entries of $M$.  For each player $i$, we use 
$\sum_{k\in I-\{i\}}d_k$ rows and $d_i$ columns of $M$.  Thus we could just use 
a totally nonsingular matrix with $D-\min_{i\in I}d_i$ rows and 
$\max_{i\in I}d_i$ columns.

For example, for a game of $3$ players with $3$ pure strategies each, using the above $6\times 6$ totally nonsingular matrix, we arrive at the system:
\begin{eqnarray*}
(\sigma_{21}+2\sigma_{22}-1)(\sigma_{31}+2\sigma_{32}-1)&=&0,\\
(2\sigma_{21}-4\sigma_{22}-1)(2\sigma_{31}-4\sigma_{32}-1)&=&0,\\
(4\sigma_{11}+16\sigma_{12}-1)(4\sigma_{31}+16\sigma_{32}-1)&=&0,\\
(8\sigma_{11}-32\sigma_{12}-1)(8\sigma_{31}-32\sigma_{32}-1)&=&0,\\
(16\sigma_{11}+128\sigma_{12}-1)(16\sigma_{21}+128\sigma_{22}-1)&=&0,\\
(32\sigma_{11}-256\sigma_{12}-1)(32\sigma_{21}-256\sigma_{22}-1)&=&0.\\
\end{eqnarray*}
If we replaced the $1$ in each factor by $\sum_{{j_k}=1}^{d_k}\sigma_{k{j_k}}$
and expanded out the polynomials, we could determine for which payoff
functions this is the system \E{_{ij}}.  For instance, in the second equation
the coefficient of $\sigma_{22}\sigma_{31}$ becomes $-9$, so this says
$u_1(s_{12},s_{22},s_{31})-u_1(s_{10},s_{22},s_{31})=-9$.  The coefficient of $\sigma_{20}\sigma_{30}$ becomes $1$, so this says 
$u_1(s_{12},s_{20},s_{30})-u_1(s_{10},s_{20},s_{30})=1$.

Now to find all the solutions to this system, we define a $D\times D$ matrix
$P$ by $P_{n(i,j)n(k,l)}=0$ if $i=k$ and $P_{n(i,j)n(k,l)}=1$ if $i\neq k$.  In
the example of a game with $3$ players, each with $3$ pure strategies, we have
\[P=
\left[\begin{array}{rrrrrr}
0&0&1&1&1&1\\
0&0&1&1&1&1\\
1&1&0&0&1&1\\
1&1&0&0&1&1\\
1&1&1&1&0&0\\
1&1&1&1&0&0\\
\end{array}
\right]\]
We associate the $n(i,j)$th row with the equation \E{_{ij}} and the $n(i,j)$th
column with the variable $\sigma_{ij}$.  So an entry of $P$ is $1$ if and only
if the corresponding variable appears in the corresponding equation.  

The reader is familiar with the \emph{determinant} of a matrix, which is the
sum of certain signed products of entries of the matrix.  The \emph{permanent}
of a matrix is the sum of those same products of entries of the matrix, but
without the signs.  In other words, for an $D\times D$ matrix $P$, the
permanent of $P$ is the sum over all permutations $\tau$ of $1,\ldots,D$ of the
products $\prod_{n=1}^DP_{n\tau(n)}$.

To find a solution of \startsystem{}, we pick $D$ entries of $P$ whose product
contributes $1$ to the permanent of $P$.  In other words, we pick a permutation
$\tau$ of $1,\ldots,D$ such that $\prod_{n=1}^DP_{n\tau(n)}=1$.  For example,
the italicized entries below represent such a choice:
\[
\left[\begin{array}{rrrrrr}
0&0&1&1&{\it1}&1\\
0&0&1&{\it1}&1&1\\
{\it1}&1&0&0&1&1\\
1&1&0&0&1&{\it1}\\
1&{\it1}&1&1&0&0\\
1&1&{\it1}&1&0&0\\
\end{array}
\right]\]
Now if $\tau(n(i,j))=n(k,l)$, this tells us to make equation 
\E{_{ij}} hold by setting the factor
$1+\sum_{{j_k}=1}^{d_k}\lambda^{(ij)}_{k{j_k}}\sigma_{k{j_k}}$ equal to zero.  
In the above example, the above choice tells us to set:
\begin{eqnarray*}
\sigma_{31}+2\sigma_{32}-1&=&0,\\
2\sigma_{21}-4\sigma_{22}-1&=&0,\\
4\sigma_{11}+16\sigma_{12}-1&=&0,\\
8\sigma_{31}-32\sigma_{32}-1&=&0,\\
16\sigma_{11}+128\sigma_{12}-1&=&0,\\
32\sigma_{21}-256\sigma_{22}-1&=&0.\\
\end{eqnarray*}
Now we have a system of $d_i$ linear equations in the $d_i$ variables
$\sigma_{ij}$, for each $i$.  Since we chose the coefficients from a totally
nonsingular matrix, each system of $d_i$ linear equations has a unique solution.
In this case, we find:
\[\sigma_{11}=\frac{7}{16},\sigma_{12}=-\frac{3}{64},
\sigma_{21}=\frac{21}{32},\sigma_{22}=\frac{5}{64},
\sigma_{31}=\frac{17}{24},\sigma_{32}=\frac{7}{48}.\]
(Clearly this particular solution does not satisfy the nonnegativity
constraints, which we would also have to check if we were interested in the
Nash equilibria of this particular game.)

Notice that this procedure would give us the same set of equations multiple
times.  For example, the choice of $D$ other entries in $P$ represented by the
italicized entries below:
\[
\left[\begin{array}{rrrrrr}
0&0&1&1&1&{\it1}\\
0&0&{\it1}&1&1&1\\
1&{\it1}&0&0&1&1\\
1&1&0&0&{\it1}&1\\
{\it1}&1&1&1&0&0\\
1&1&1&{\it1}&0&0\\
\end{array}
\right]\]
gives the same system of equations.  The problem is that for each set of $d_i$
columns corresponding to the variables $\sigma_{i1},\ldots,\sigma_{i{d_i}}$,
we can apply any permutation to that set of columns without affecting the
meaning of our choice.  So if we carried out this procedure na\"\i vely, it
would repeat each solution $\prod_{i\in I}d_i!$ times.  We should avoid solving
the same system twice.  However, if we obtain the same solution from a
different choice of which factor in each of the equations to set to zero, then
we should perturb our totally nonsingular matrix so this doesn't happen, for
reasons which will become clear later.

Carrying out this procedure, we find all $10$ roots of \E{}.  We list the
values of
$(\sigma_{11},\sigma_{12},\sigma_{21},\sigma_{22},\sigma_{31},\sigma_{32})$
below, along with the corresponding permutations of $1,\ldots,D$:
\begin{eqnarray*}
\left(\frac{3}{64},\frac{1}{512},\frac{3}{4},
\frac{1}{8},\frac{3}{16},\frac{1}{64}\right)&\qquad&5,6,1,2,3,4\\
\left(\frac{7}{32},\frac{3}{128},\frac{21}{16},
\frac{-5}{32},\frac{5}{12},\frac{-1}{24}\right),&\qquad&4,6,1,5,2,3\\
\left(\frac{17}{96},\frac{7}{384},\frac{21}{16},
\frac{-5}{32},\frac{7}{8},\frac{3}{16}\right)&\qquad&3,6,1,5,2,4\\
\left(\frac{5}{48},\frac{-1}{192},\frac{129}{160},
\frac{31}{320},\frac{5}{12},\frac{-1}{24}\right)&\qquad&4,5,1,6,2,3\\
\left(\frac{7}{16},\frac{-3}{64},\frac{129}{160},
\frac{31}{320},\frac{7}{8},\frac{3}{16}\right)&\qquad&3,5,1,6,2,4\\
\left(\frac{7}{32},\frac{3}{128},\frac{33}{80},
\frac{-7}{160},\frac{7}{4},\frac{-3}{8}\right)&\qquad&4,6,2,5,1,3\\
\left(\frac{17}{96},\frac{7}{384},\frac{33}{80},
\frac{-7}{160},\frac{17}{24},\frac{7}{48}\right)&\qquad&3,6,2,5,1,4\\
\left(\frac{5}{48},\frac{-1}{192},\frac{21}{32},
\frac{5}{64},\frac{7}{4},\frac{-3}{8}\right)&\qquad&4,5,2,6,1,3\\
\left(\frac{7}{16},\frac{-3}{64},\frac{21}{32},
\frac{5}{64},\frac{17}{24},\frac{7}{48}\right)&\qquad&3,5,2,6,1,4\\
\left(\frac{3}{16},\frac{1}{64},\frac{3}{64},
\frac{1}{512},\frac{3}{4},\frac{1}{8}\right)&\qquad&3,4,5,6,1,2
\end{eqnarray*}
We note that the first and last of these satisfy the nonnegativity constraints,
and hence are Nash equilibria.  (The reader may notice some symmetries between
pairs of the solutions.  This is because we happened to use the same region of
$M$ for all the players.)

We can find the rest of the solutions of the system $(*)$ corresponding to this
same game too.  Suppose we require $v_{ij}$ to be positive instead of
vanishing.  Then the variable $\sigma_{ij}$ goes away from every equation and
the equation \E{_{ij}} is replaced by the equation $\sigma_{ij}=0$.  In the
above example, suppose we require $v_{31}$ to be positive.  Then we obtain a
new system
\begin{eqnarray*}
(\sigma_{21}+2\sigma_{22}-1)(2\sigma_{32}-1)&=&0,\\
(2\sigma_{21}-4\sigma_{22}-1)(4\sigma_{32}-1)&=&0,\\
(4\sigma_{11}+16\sigma_{12}-1)(16\sigma_{32}-1)&=&0,\\
(8\sigma_{11}-32\sigma_{12}-1)(32\sigma_{32}-1)&=&0,\\
\sigma_{31}&=&0,\\
(32\sigma_{11}-256\sigma_{12}-1)(32\sigma_{21}-256\sigma_{22}-1)&=&0.\\
\end{eqnarray*}
This system corresponds to the $n(i,j),n(i,j)$ minor of $M$, which is also
totally nonsingular.  We can use the same matrix $P$, but zero out the
$n(i,j)$th row and $n(i,j)$th column and set the $n(i,j),n(i,j)$ entry to 1.  
In our example:
\[P'_{ij}=
\left[\begin{array}{rrrrrr}
0&0&1&1&0&1\\
0&0&1&1&0&1\\
1&1&0&0&0&1\\
1&1&0&0&0&1\\
0&0&0&0&1&0\\
1&1&1&1&0&0\\
\end{array}
\right]\]
Now we follow the same procedure as above, performing the cofactor expansion of
the permanent along the $n(i,j)$th row.  For example, the product of italicized
entries below contributes $1$ to the permanent of this matrix:
\[
\left[\begin{array}{rrrrrr}
0&0&1&{\it1}&0&1\\
0&0&{\it1}&1&0&1\\
{\it1}&1&0&0&0&1\\
1&1&0&0&0&{\it1}\\
0&0&0&0&{\it1}&0\\
1&{\it1}&1&1&0&0\\
\end{array}
\right]\]
The corresponding system is:
\begin{eqnarray*}
\sigma_{21}+2\sigma_{22}-1&=&0,\\
2\sigma_{21}-4\sigma_{22}-1&=&0,\\
4\sigma_{11}+16\sigma_{12}-1&=&0,\\
32\sigma_{32}-1&=&0,\\
\sigma_{31}&=&0,\\
32\sigma_{11}-256\sigma_{12}-1&=&0.\\
\end{eqnarray*}
Its solution is:
\begin{eqnarray*}
\sigma_{11}=\frac{17}{96},\sigma_{12}=\frac{7}{384},\\
\sigma_{21}=\frac{3}{4},\sigma_{22}=\frac{1}{8},\\
\sigma_{31}=0,\sigma_{32}=\frac{1}{32}.\\
\end{eqnarray*}
This at least satisfies the nonnegativity constraints on the $\sigma_{ij}$'s.
If we are interested in the Nash equilibria of this game,
we also have to check that it satisfies the nonnegativity constraints on
the $v_{ij}$'s, namely, that $v_{31}\geq0$.  We substitute the $\sigma_{ij}$'s
into
$$(16\sigma_{11}+128\sigma_{12}-1)(16\sigma_{21}+128\sigma_{22}-1),$$
the expected payoff to player $3$ from playing $\sigma_{31}$,
obtaining $\frac{225}{2}$, which is strictly greater than zero.  So this
solution to the polynomial system is \emph{not} a Nash equilibrium, since
$s_{31}$ is a strictly better response to $\sigma_{-3}$ than the value of
$\sigma_3$ given by this solution.

In this way we see that our specially constructed factorizable game of a given
format $(N;d_1,\ldots,d_N)$ contains subgames of every smaller format
$(N':d'_1,\ldots,d'_{N'})$, with $N'\leq N$ and $d'_i\leq d_i$ for each $i$,
such that the subgames are also factorizable.  As we will see, we only need to
solve one (or a few) of the factorizable polynomial systems for each format in
the manner we have described so far.

Geometrically, by constructing the $D\times D$ totally nonsingular matrix $M$,
we picked $D$ vectors in $D$-dimensional space, 
such that not only are these all distinct points, but if we project any $m$ of
them onto any $m$-dimensional coordinate subspace, the images are all also
distinct.  The condition that two particular such images coincide is an
equation, which is satisfied only on a subset of real $D$-dimensional space of
strictly lower dimension.  So every open subset of real $D$-dimensional space
does \emph{not} satisfy the condition, almost everywhere.  Since there are only
finitely many of these conditions, every open subset does not satisfy any of
them, almost everywhere.  In particular, we could construct a totally
nonsingular matrix such that all the nonnegativity constraints of our
factorizable game also held with strict inequality.  However,
while this would provide an example of a game with the maximal possible number
of totally mixed Nash equilibria, it would not be particularly
relevant to the use we will be making of our specially constructed factorizable
games.

\section{Polyhedra and Polynomial Systems}
\noindent
With a system of polynomial equations is associated a \emph{polyhedral
subdivision}, that is, a polyhedron which is subdivided into \emph{cells}, each
of which is also a polyhedron, glued together along their faces.  We illustrate
this for a game of 3 players with 2 pure strategies each, since in this case
the polyhedral subdivision is 3-dimensional.  Using the same totally
nonsingular matrix, we obtain the following system of factored equations:
\begin{eqnarray}
\label{factored23}
\color{green}(\sigma_{21}-1)(\sigma_{31}-1)&\color{green}=&\color{green}0,\\
\label{factored13}
\color{red}(2\sigma_{11}-1)(2\sigma_{31}-1)&\color{red}=&\color{red}0,\\
\label{factored12}
\color{blue}(4\sigma_{11}-1)(4\sigma_{21}-1)&\color{blue}=&\color{blue}0.
\end{eqnarray}
Expanding this out, we obtain
\begin{eqnarray}
\label{expanded23}
\color{green}\sigma_{21}\sigma_{31}-\sigma_{21}-\sigma_{31}+1
&\color{green}=&\color{green}0,\\
\label{expanded13}
\color{red}4\sigma_{11}\sigma_{31}-2\sigma_{11}-2\sigma_{31}+1&
\color{red}=&\color{red}0,\\
\label{expanded12}
\color{blue}16\sigma_{11}\sigma_{21}-4\sigma_{11}-4\sigma_{21}+1&
\color{blue}=&\color{blue}0.
\end{eqnarray}

\begin{figure}
\begin{center}
\includegraphics{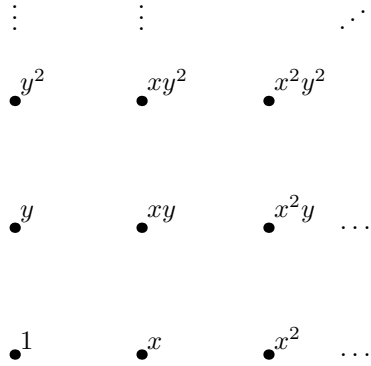}
\end{center}
\caption{The Monomial Lattice in Two Variables}
\label{monomial_lattice}
\end{figure}

A monomial $x_1^{\alpha_1}x_2^{\alpha_2}\cdots x_n^{\alpha_n}$ in $n$ variables
can be represented by the \emph{lattice point}
$(\alpha_1,\ldots,\alpha_n)\in\N^n$ of its exponents.  For example, the lattice
of monomials in two variables $x$ and $y$ is depicted in Figure
\ref{monomial_lattice}.

\begin{figure}
\begin{center}
\includegraphics{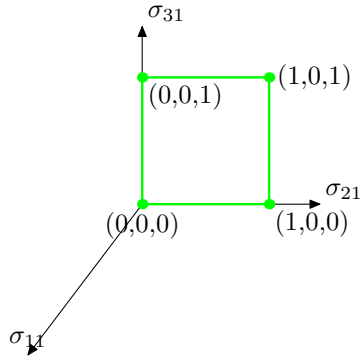}
\end{center}
\caption{The Newton Polytope of 
$\sigma_{21}\sigma_{31}-\sigma_{21}-\sigma_{31}+1=0$}
\label{newton23}
\end{figure}

\begin{figure}
\begin{center}
\includegraphics{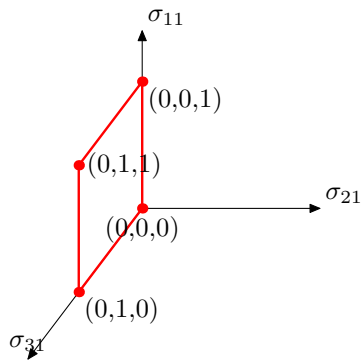}
\end{center}
\caption{The Newton Polytope of 
$4\sigma_{11}\sigma_{31}-2\sigma_{11}-2\sigma_{31}+1=0$}
\label{newton13}
\end{figure}

\begin{figure}
\begin{center}
\includegraphics{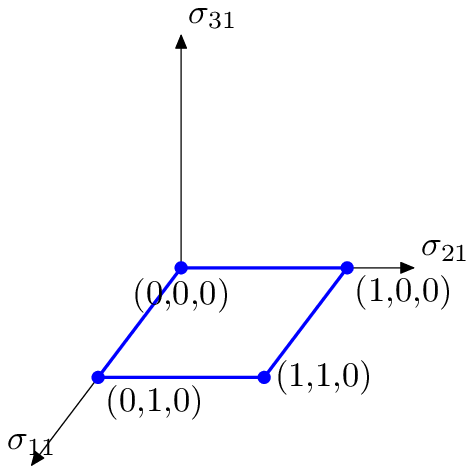}
\end{center}
\caption{The Newton Polytope of 
$16\sigma_{11}\sigma_{31}-4\sigma_{11}-4\sigma_{31}+1=0$}
\label{newton12}
\end{figure}

The \emph{Newton
polytope} of a polynomial equation is the convex hull of the lattice points of
the monomials occuring in that equation.  In our example system in 3 variables
$\sigma_{11}$, $\sigma_{21}$, and $\sigma_{31}$, the Newton polytope of
Equation \ref{expanded23} is depicted in Figure \ref{newton23},  
the Newton polytope of
Equation \ref{expanded13} is depicted in Figure \ref{newton13},
and the Newton
polytope of Equation \ref{expanded12} is depicted in Figure \ref{newton12}.

\begin{figure}
\begin{center}
\includegraphics{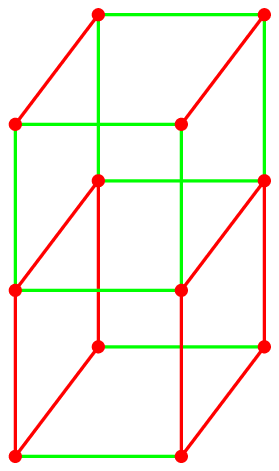}
\end{center}
\caption{Minkowski Sum of A Pair Of Newton Polytopes}
\label{minkowskisum}
\end{figure}

The \emph{Minkowski sum} of an $n$-dimensional polytope with vertices
$V_{01},\ldots,V_{0{m_0}}$ and an $n$-dimensional polytope with vertices 
$V_{11},\ldots,V_{1{m_1}}$ is the convex hull of the points $V_{0i}+V_{1j}$ in
$n$-dimensional space, for $i=1,\ldots,m_0$ and $j=1,\ldots,m_1$.  Figure
\ref{minkowskisum} depicts the Minkowski sum of the Newton polytopes of Equation
\ref{expanded23} and \ref{expanded13}.
We can think of the Minkowski sum as translating one of the polytopes along
each edge of the other polytope.  Here we have colored the vertices and edges
as if we first translated the red polytope along each edge of the green
polytope.  The red edges came from the original red polytope, and the green
edges came from edges of the green polytope along which we translated.  Notice
that we can do this in more than one way.  For example, we could have colored
the vertices and edges as if we translated the green polytope along each edge
of the red polytope first.  Such a coloring of the Minkowski sum gives us a
\emph{mixed subdivision}, which in this case has two cells, the two cubes
in Figure \ref{minkowskisum}.

\begin{figure}
\begin{center}
\includegraphics{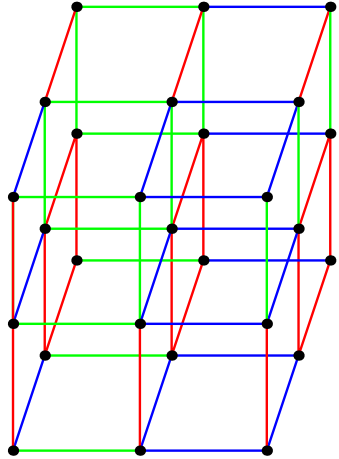}
\end{center}
\caption{Minkowski Sum of Three Newton Polytopes}
\label{mixedvolume}
\end{figure}

Finally, the Minkowski sum of all three of our Newton polytopes is depicted in
Figure \ref{mixedvolume}.
A cell of a mixed subdivision is \emph{mixed} if
each color delineates only edges (or possibly vertices) in the cell,
not higher-dimensional faces.  In Figure \ref{mixedvolume}, the top left front
cube is not mixed, because two of its faces are green squares, and the top
right front cube is not mixed, because two of its faces are blue squares.
We see that two of the cells in this mixed subdivision are mixed.  Each mixed
cell tells us how to obtain certain solutions to the factorizable polynomial
system.  Namely, in each polynomial equation, we should look at the edges with
the corresponding color in the mixed cell, and set the factor(s) corresponding
to the directions of those edges to zero.  In the game-theoretic case, there
will be exactly one solution corresponding to each mixed cell, since the
polynomial system has degree at most $1$ in any variable.

\begin{figure}
\begin{center}
\includegraphics{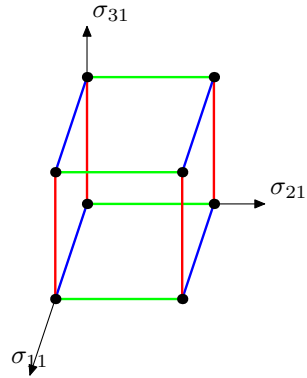}
\end{center}
\caption{One Mixed Cell of the Mixed Subdivision}
\label{bottommixedcell}
\end{figure}

For example, Figure \ref{bottommixedcell} depicts the bottom mixed cell
of the mixed subdivision.
To find the solution corresponding to this mixed cell, we should set
the factor {\color{green} $\sigma_{21}-1$} to zero in the green 
Equation \ref{factored23}; we
should set the factor {\color{red} $2\sigma_{31}-1$} to zero in the red 
Equation \ref{factored13}; and we should set the factor 
{\color{blue} $4\sigma_{11}-1$} to zero in the
blue Equation \ref{factored12}.  This gives us the solution
\[{\color{blue}\sigma_{11}=\frac{1}{4}},
{\color{green}\sigma_{21}=1},
{\color{red}\sigma_{31}=\frac{1}{2}}.\]

\begin{figure}
\begin{center}
\includegraphics{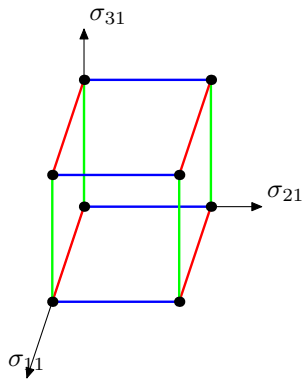}
\end{center}
\caption{Another Mixed Cell of the Mixed Subdivision}
\label{topmixedcell}
\end{figure}

Figure \ref{topmixedcell} depicts the top mixed cell of the mixed
subdivision.
To find the solution corresponding to this mixed cell, we should set
the factor {\color{green}$\sigma_{31}-1$} to zero in the green Equation
\ref{factored23}; we
should set the factor {\color{red}$2\sigma_{11}-1$} to zero in the red Equation
\ref{factored13}; and we should set the factor {\color{blue}$4\sigma_{21}-1$} 
to zero in the blue Equation \ref{factored12}.  This gives us the solution
\[{\color{red}\sigma_{11}=\frac{1}{2}},
{\color{blue}\sigma_{21}=\frac{1}{4}},
{\color{green}\sigma_{31}=1}.\]
\noindent
In the game-theoretic case, the cells of the mixed subdivision are always 
$D$-dimensional cubes (or hypercubes).  The matrix $P$ in this example is
\[\left[\begin{array}{rrr}
\color{green}0&\color{green}1&\color{green}1\\
\color{red}1&\color{red}0&\color{red}1\\
\color{blue}1&\color{blue}1&\color{blue}0
\end{array}\right]\]
Recall that the rows correspond to equations (and hence each will have a
different color), and the columns correspond to variables.
A permutation which contributes to the permanent will tell us how to color each
variable.  The edges of the mixed cube going in the direction corresponding to
that variable will be colored with that color, and the linear factor containing
that variable in that equation will be set to zero.

\section{Finding All Nash Equilibria by Polyhedral Homotopy Continuation}
\noindent
We have taken considerable trouble to find all Nash equilibria of one
particular kind of specially constructed game.  The good news is that once we
have done this for one particular game format, we can easily solve \emph{any}
generic game of that format.\footnote{It's not quite true that we only need to
solve \emph{one} game of each format.  Something untoward could happen on the
way from this game to the one we are interested in; e.g., we could run into a
game whose set of Nash equilibria is positive-dimensional.  So we should have
the solutions to a few of the specially constructed games ready to hand.  We
can make more such games by using different regions of our totally nonsingular
matrix.} Furthermore, we can even more easily look for one or a few of the Nash
equilibria, look for Nash equilibria with some particular small support, and so
forth.  The key idea is to ``morph'' the specially constructed polynomial
system into the polynomial system we are actually interested in.  As we do
this, the solutions to \E{} for the specially constructed game will also morph
into solutions to \E{} for the
game of interest.  Moreover, this procedure is ``embarrassingly parallel''.
The morphing of each solution is independent of the morphing of every other
solution.  We can partition the solutions to \E{} for the specially constructed
game into subsets, and hand each subset to a different processor.  If somewhere
along the way a morphed solution begins to look disappointing (for example, it
doesn't look like it will end up satisfying the nonnegativity constraints, or
it starts to have an imaginary component that we fear won't go away), we can
always stop morphing that solution and come back to it later if more promising
ones don't pan out.

Readers may be familiar with the homotopy continuation method under the guise
of ``tracing procedures'', such as the techniques of Lemke-Howson or Govindan
and Wilson.  To solve a polynomial system by the homotopy continuation method,
we create a family $\mathcal F$ of polynomial systems
$f_{1t}=0,f_{2t}=0,\ldots,f_{mt}=0$ parametrized by a variable $t$ lying in
$[0,1]$, such that the polynomial system we want to solve is
$f_{11}=0,f_{21}=0,\ldots,f_{m1}=0$, and the polynomial system
$f_{10}=0,f_{20}=0,\ldots,f_{m0}=0$, called the \emph{starting system}, is easy
to solve.  We consider each polynomial system in $\mathcal F$ to lie in some
topological space such that for any particular point, the map of the polynomial
space which evaluates the polynomial at that point is continuous.  We require
the map from $[0,1]$ to this space to be continuous, or in other words a
\emph{homotopy}.  Now suppose
$(x_1,\ldots,x_n)$ satisfies the polynomial system 
$f_{1{t_0}},\ldots,f_{m{t_0}}$, and $t_1$ is near $t_0$.  Since the homotopy is
continuous, $f_{1{t_1}},\ldots,f_{m{t_1}}$ must be near
$f_{1{t_0}},\ldots,f_{m{t_0}}$, and so $f_{1{t_1}}(x_1,\ldots,x_n)$ must lie
near zero, $f_{2{t_1}}(x_1,\ldots,x_n)$ must lie near zero, and so forth.
Therefore, since polynomial functions are also continuous, 
we can look for a root of the system $f_{1{t_1}},\ldots,f_{m{t_1}}$
near $(x_1,\ldots,x_n)$.  We make a \emph{prediction}, i.e., we guess a
possible root of $f_{1{t_1}},\ldots,f_{m{t_1}}$ near $(x_1,\ldots,x_n)$, and
then a \emph{correction}, i.e., we find an actual root near our guess, using
Newton's method for example.  Once we have one, we can proceed to the next
iteration for $t_2$ near $t_1$, and so forth.  At the end we will have a path
from our original root $(x_1,\ldots,x_n)_{t=0}$ to a root
$(x_1,\ldots,x_n)_{t=1}$ of the desired system.  \cite{heringspeetershomotopy}
presents a survey of previous uses of homotopy continuation methods in game
theory.  The book \cite{SommeseWampler} gives a recent survey of numerical
methods for solving polynomial systems, including detailed treatment of
homotopy continuation, and in particular polyhedral homotopy continuation.

If we fix the number of equations, and the Newton polytopes of each equation,
then the set of such polynomial systems becomes a vector space over the
coefficient field.  Each monomial occurring in each equation corresponds to a
basis element of this vector space, and a particular polynomial system is
uniquely specified by giving the coefficients of all the monomials in all the
equations.  In particular, if the coefficient field is $\R$ then this space of
polynomial systems is a finite-dimensional real vector space.  Hence,
it is equipped with a topology, the usual topology of such spaces.  We will
call the number of equations together with the Newton polytopes of each
equation the \emph{shape} of a polynomial system.  Polyhedral homotopy
continuation is simply homotopy continuation among polynomial systems of the
same shape.  The word ``polyhedral'' refers to the polyhedral subdivision
introduced in the previous section.  

The Bernstein-Kouchnirenko theorem \cite{bernstein} \cite{Kouchnirenko} tells us
that the number of $0$-dimensional complex roots, none of whose components are
zero, of every generic polynomial system of a given shape is the same.  This
number is called the \emph{Bernstein number} of the system.  Polyhedral
homotopy continuation \cite{hubersturmfels} provides an alternate constructive
proof of this fact.  Thus, if we apply polyhedral homotopy continuation to a 
generic polynomial system we
will find a unique (possibly complex) root of the system in question at the end
of each path leading from one of the roots of the starting system.
\cite{mckelveymclennan} applied the Bernstein-Kouchnirenko theorem to find the
number of complex roots of the polynomial system \E{} for a generic game: it is
the permanent of the matrix $P$, divided by $\prod_{i\in I}d_i!$.  In
\cite{dattathesis} we generalized this theorem to polynomial systems obeying
special conditions, such as those arising from graphical games.  Briefly put,
the special conditions imply that some of the entries in $P$ are zero.  

In general, the most difficult part of polyhedral homotopy continuation is
computing the mixed subdivision, finding a polynomial system which is generic,
and using the mixed subdivision to find all the roots.  Once all this is done,
we can follow a simple linear homotopy (i.e., of the form
$f_t=(1-t^k)f_0+t^kf_1$ for some $k\geq 1$)
from this starting system to the desired
system, which is relatively straightforward.  Therefore polyhedral homotopy
continuation is particularly well-suited in the case of game theory, since we
know exactly how to find and solve a starting system.  The system \E{_{ij}} is
an example of a \emph{multihomogeneous} system.  Every monomial appearing in
one of these equations has the same degree, namely $1$, in all of player $1$'s
variables put together, the same degree, namely $1$, in in all of player $2$'s
variables put together, and so forth.  In the same way multihomogeneous systems
are generally easy to solve by polyhedral homotopy continuation.  As we saw,
adding in the conditions \E{_{i0}} does not make the problem more difficult,
although it does mean multihomogeneity no longer holds.  The system \E{} is a
\emph{linear product family}, as described in Section 8.4.3 of
\cite{SommeseWampler}.

We may not even need to find all the roots of the starting system before
starting to find the Nash equilibria of the desired game.  Once we have a
single root of the starting system, we can start tracing it.  However, if we
are unlucky the corresponding root of the desired game may not be nonnegative
or even real, in which case we will have to go find another root of the
starting system to trace.  \cite{McLennanRealRoots} describes the expected
number of real roots of a random multihomogeneous system of polynomial
equations, and gives the lower bound as the square root of the Bernstein number
(the number of complex roots of a generic game).  

Let's expand out our factorizable polynomial system for the game of 3 players
with 3 pure strategies each.
\begin{eqnarray*}
\sigma_{21}\sigma_{31}+2\sigma_{21}\sigma_{32}+2\sigma_{22}\sigma_{31}
+4\sigma_{22}\sigma_{32}\\
\qquad-\sigma_{21}-2\sigma_{22}-\sigma_{31}-2\sigma_{32}+1
&=&0,\\
4\sigma_{21}\sigma_{31}-8\sigma_{21}\sigma_{32}-8\sigma_{22}\sigma_{31}
+16\sigma_{22}\sigma_{32}\\
\qquad-2\sigma_{21}+4\sigma_{22}-2\sigma_{31}+4\sigma_{32}+1
&=&0,\\
16\sigma_{11}\sigma_{31}+64\sigma_{11}\sigma_{32}+64\sigma_{12}\sigma_{31}
+256\sigma_{12}\sigma_{32}\\
\qquad-4\sigma_{11}-16\sigma_{12}-4\sigma_{31}-16\sigma_{32}
+1&=&0,\\
64\sigma_{11}\sigma_{31}+256\sigma_{11}\sigma_{32}+256\sigma_{12}\sigma_{31}
+1024\sigma_{12}\sigma_{32}\\
\qquad-8\sigma_{11}-32\sigma_{12}-8\sigma_{31}-32\sigma_{32}+1&=&0,\\
256\sigma_{11}\sigma_{21}+2048\sigma_{11}\sigma_{22}+2048\sigma_{12}\sigma_{21}
+16384\sigma_{12}\sigma_{22}\\
\qquad-16\sigma_{11}-128\sigma_{12}-16\sigma_{21}-128\sigma_{22}+1&=&0,\\
1024\sigma_{11}\sigma_{21}-8192\sigma_{11}\sigma_{22}-8192\sigma_{12}\sigma_{21}
+65536\sigma_{12}\sigma_{22}\\
\qquad-32\sigma_{11}+256\sigma_{12}-32\sigma_{21}+256\sigma_{22}+1&=&0.
\end{eqnarray*}
We can make an input file for the polyhedral homotopy continuation software
\PHC{pack} \cite{phc} which specifies this polynomial system.  Since \PHC{pack}
orders the variables according to the order they have appeared in the file, we
will write the equations in the reverse order from the system listed above, so
that the variables will appear in the order
$\sigma_{11},\sigma_{12},\sigma_{21},\sigma_{22},\sigma_{31},\sigma_{32}$.  
Here is the input file 
{\tt gameof3x3x3\_start\_phc}:
\begin{verbatim}
6
1 - 32*s11 + 256*s12 - 32*s21 + 256*s22
+ 1024*s11*s21 - 8192*s11*s22 - 8192*s12*s21 + 65536*s12*s22;
1 - 16*s11 - 128*s12 - 16*s21 - 128*s22 
+ 256*s11*s21 + 2048*s11*s22 + 2048*s12*s21 + 16384*s12*s22;
1 - 8*s11 - 32*s12 - 8*s31 - 32*s32
+ 64*s11*s31 + 256*s11*s32 + 256*s12*s31 + 1024*s12*s32;
1 - 4*s11 - 16*s12 - 4*s31 - 16*s32
+ 16*s11*s31 + 64*s11*s32 + 64*s12*s31 + 256*s12*s32;
1 - 2*s21 + 4*s22 - 2*s31 + 4*s32
+ 4*s21*s31 - 8*s21*s32 - 8*s22*s31 + 16*s22*s32;
1 - s21 - 2*s22 - s31 - 2*s32
+ s21*s31 + 2*s21*s32 + 2*s22*s31 + 4*s22*s32;
\end{verbatim}
The first line specifies the number of equations, and the rest of the file
specifies the equations.  Since an unknown can only consist of up to 5
characters, we denote the variable $\sigma_{11}$ by {\tt s11} in the file, and
so forth.

We could find the roots of this factorizable system using
\PHC{pack} itself (in which case we should tell \PHC{} about the linear product
structure of our equations), or using other programs in the manner described
above.  We will list the 10 roots of the starting system we found before in
another input file for \PHC{pack}, which we call {\tt
gameof3x3x3\_start\_phc.roots}.  Here is the beginning of that file, including
the first two roots:
\begin{verbatim}
10 6
===========================================================
solution 1 :
t :  0.00000000000000E+00   0.00000000000000E+00
m : 1
the solution for t :
 s11 :  4.68750000000000e-02   0.00000000000000E+00
 s12 :  1.95312500000000e-03   0.00000000000000E+00
 s21 :  7.50000000000000e-01   0.00000000000000E+00
 s22 :  1.25000000000000e-01   0.00000000000000E+00
 s31 :  1.87500000000000e-01   0.00000000000000E+00
 s32 :  1.56250000000000e-02   0.00000000000000E+00
== err :  0.000E+00 = rco :  1.000E+00 = res :  0.000E+00 ==
solution 2 :
t :  0.00000000000000E+00   0.00000000000000E+00
m : 1
the solution for t :
 s11 :  2.18750000000000e-01   0.00000000000000E+00
 s12 :  2.34375000000000e-02   0.00000000000000E+00
 s21 :  1.31250000000000e+00   0.00000000000000E+00
 s22 : -1.56250000000000e-01   0.00000000000000E+00
 s31 :  4.16666666666667e-01   0.00000000000000E+00
 s32 : -4.16666666666667e-02   0.00000000000000E+00
== err :  0.000E+00 = rco :  1.000E+00 = res :  0.000E+00 ==
\end{verbatim}
The first line indicates that the file contains $10$ solutions in $6$ unknowns.
Here $t$ denotes the homotopy parameter, and $m$ denotes the multiplicity of
each root.  A line such as
\begin{verbatim}
 s11 :  4.68750000000000e-02   0.00000000000000E+00
\end{verbatim}
indicates that at this solution, the variable {\tt s11} has real part $.046875$
and imaginary part $0$.  The lines
\begin{verbatim}
the solution for t :
\end{verbatim}
and
\begin{verbatim}
== err :  0.000E+00 = rco :  1.000E+00 = res :  0.000E+00 ==
\end{verbatim}
are lines that would have been included by \PHC{pack} if it had written this
solution file itself, so we include them also even if we didn't use \PHC{pack}
to generate these solutions.

Finally, we write an input file {\tt gameof3x3x3\_phc} with an example of
another polynomial system of the same shape, which we would like to solve:
\begin{verbatim}
6
1 - 2*s11 + 3*s12 - 5*s21 + 7*s22
- 7*s11*s21 - 5*s11*s22 - 3*s12*s21 + 2*s12*s22;
7 - 3*s11 - 5*s12 + 2*s21 - 3*s22 
- 7*s11*s21 + 3*s11*s22 + s12*s21 - s12*s22;
3 - 5*s11 - 3*s12 - 2*s31 + 2*s32
+ 5*s11*s31 + 7*s11*s32 - 7*s12*s31 + s12*s32;
2 - 3*s11 - 5*s12 - 7*s31 + 7*s32
+ 5*s11*s31 + 3*s11*s32 - 2*s12*s31 - s12*s32;
1 - 2*s21 - 3*s22 + 7*s31 - 5*s32
- s21*s31 + 2*s21*s32 + 5*s22*s31 + 3*s22*s32;
1 - s21 + 2*s22 - 3*s31 - 5*s32
+ 7*s21*s31 - 2*s21*s32 + 5*s22*s31 + 3*s22*s32;
\end{verbatim}

Now we can invoke \PHC{pack} with the {\tt -p} option, indicating that we
already have a starting system and its solutions.
\begin{verbatim}
$ phc -p
Welcome to PHC (Polynomial Homotopy Continuation) V2.3.16 25 Nov 2006
Polynomial Continuation defined by a homotopy in one parameter.

Reading the target polynomial system...
Give a string of characters : gameof3x3x3_phc

Reading the name of the output file.
Give a string of characters : gameof3x3x3_phc.output

Do you want the solutions on separate file ? (y/n) y
Reading the name of the file to write the solutions on.
Give a string of characters : gameof3x3x3_phc.roots

Reading the name of the file for start system.
Give a string of characters : gameof3x3x3_start_phc

Reading the name of the file for the solutions.
Give a string of characters : gameof3x3x3_start_phc.roots

Homotopy is H(x,t) = a*(1-t)^k * Q(x) + t^k * P(x) = 0, t in [0,1],
      with Q(x) = 0 a start system, and P(x) = 0 the target system.
\end{verbatim}
At this point we are presented with several menus allowing us to change
different options for controlling the homotopy continuation.  At each point we
can enter {\tt 0} to accept the default options.  Finally we are presented
with:
\begin{verbatim}
No more input expected.  See output file for results.
\end{verbatim}
and, possibly after some delay depending on how big our system is (in this
case, there is no noticeable delay), the program exits.  Now we can look at the
solution file we specified {\tt gameof3x3x3\_phc.roots}:
\begin{verbatim}
10 6
===========================================================
solution 1 :
t :  1.00000000000000E+00   0.00000000000000E+00
m : 1
the solution for t :
 s11 :  5.52632039981343E-01   1.29171976073713E+00
 s12 : -1.53240191264371E+00  -4.25943189689837E-01
 s21 :  9.41581793200150E-02  -6.89962841355063E-01
 s22 :  4.25278257933885E-03   1.43036622569914E+00
 s31 :  2.32439674939101E-01   4.46695512464090E-01
 s32 : -3.51970339164687E-01   2.04469018836310E-01
== err :  8.996E-16 = rco :  7.083E-02 = res :  7.383E-15 ==
\end{verbatim}
We see that each component of this solution has a nonzero imaginary part.  So
this solution is not of interest to us.  Looking further down in the file, we
see another solution:
\begin{verbatim}
solution 3 :
t :  1.00000000000000E+00   0.00000000000000E+00
m : 1
the solution for t :
 s11 :  1.27522488578381E+00   0.00000000000000E+00
 s12 :  7.45738698011832E-01  -3.26265223399926E-55
 s21 : -1.04186142941727E-01   4.07831529249908E-55
 s22 : -1.12076297688423E+00   6.52530446799852E-55
 s31 : -5.09803187724616E-01  -1.02304887506437E-55
 s32 :  4.44045922481355E-01  -2.65090494012440E-55
== err :  5.009E-16 = rco :  6.629E-02 = res :  3.664E-15 ==
\end{verbatim}
Here the imaginary parts occurring in the various components are very small,
and could be due to numerical error.  The eighth solution is similar.  To test
our hypothesis, we make another file {\tt gameof3x3x3\_phc.real\_roots} in
which we include only these two roots, setting their imaginary parts to zero
and renumbering them in sequence:
\begin{verbatim}
2 6
===========================================================
solution 1 :
t :  1.00000000000000E+00   0.00000000000000E+00
m : 1
the solution for t :
 s11 :  1.27522488578381E+00   0.00000000000000E+00
 s12 :  7.45738698011832E-01   0.00000000000000E+00
 s21 : -1.04186142941727E-01   0.00000000000000E+00
 s22 : -1.12076297688423E+00   0.00000000000000E+00
 s31 : -5.09803187724616E-01   0.00000000000000E+00
 s32 :  4.44045922481355E-01   0.00000000000000E+00
== err :  5.009E-16 = rco :  6.629E-02 = res :  3.664E-15 ==
solution 2 :
t :  1.00000000000000E+00   0.00000000000000E+00
m : 1
the solution for t :
 s11 :  6.39293179706243E-02   0.00000000000000E+00
 s12 : -2.16568143357771E+00   0.00000000000000E+00
 s21 :  4.93650795841189E+01   0.00000000000000E+00
 s22 : -1.96619254862997E+01   0.00000000000000E+00
 s31 : -6.49203588219902E-01   0.00000000000000E+00
 s32 : -1.51339980038990E+00   0.00000000000000E+00
== err :  2.780E-13 = rco :  1.820E-05 = res :  1.670E-13 ==
\end{verbatim}
Then we ask \PHC{pack} to validate them by calling it with the {\tt -v} option:
\begin{verbatim}
$ phc -v
Welcome to PHC (Polynomial Homotopy Continuation) V2.3.16 25 Nov 2006
Validation, refinement and purification of computed solution lists.

MENU with Validation Methods : 
  0. Scanning (huge) solution files and creating condition tables;
  1. Basic Validation : refining and weeding out the solution set;
  2. Evaluation of the residuals using multi-precision arithmetic;
  3. Newton's method using multi-precision arithmetic;
  4. Winding-Number Computation by homotopy continuation;
  5. Polyhedral Validation : frequency table of path directions;
  6. Newton's method with deflation for isolated singularities;
  7. Multiplicity structure of isolated singular solutions.
Type 0, 1, 2, 3, 4, 5, 6, or 7 to select, or i for info : 2

Is the system on a file ? (y/n/i=info) y

Reading the name of the input file.
Give a string of characters : gameof3x3x3_phc

Reading the name of the output file.
Give a string of characters : gameof3x3x3_phc.validation_of_real_roots

Reading the name of the file for the solutions.
Give a string of characters : gameof3x3x3_phc.real_roots

Give the number of decimal places : 16

\end{verbatim}
In this case we chose to evaluate each polynomial system at our candidate
roots, obtaining the \emph{residuals} (the magnitudes of their images, which
were supposed to vanish).
We look in the file {\tt gameof3x3x3\_phc.validation\_of\_real\_roots} for the
section beginning {\tt THE RESIDUALS}:
\begin{verbatim}
THE RESIDUALS with 16 decimal places :
residual 1 : 1.2838672747E-14
residual 2 : 6.530598545E-13
\end{verbatim}
So we suspect that these roots are indeed real.  Since they do not satisfy the
nonnegativity constraints, they are not Nash equilibria.  (Note well that even
if all the components of a given solution were real and nonnegative, we would
still have to check that {\tt s11 + s12 <= 1}, {\tt s21 + s22 <= 1}, and {\tt
s31 + s32 <= 1}.)

Thus, we can compute a library of starting systems for games of various
formats.  This initial computation may take a long time (indeed, computing the
permanent is $NP$-hard), but only has to be done once
(or a few times) for games of each format.  Once we have done this, for any
given game, we can look in the library for an appropriate starting system along
with its roots, parcel out the roots among the processors we are using
(possibly according to some heuristic scheme if we believe some of them are
more likely to lead to Nash equilibria), and hand the starting system, the
desired system, and the subset of roots to each processor.

If we wish, we can use interval computation along the path to each root to get
a verified bound on where the final root is.  However, interval techniques
become computationally expensive in higher dimensions.  For example, in $D$
dimensions an ``interval'' might be a box with $2^D$ corners.  

\section{Solving Polynomial Systems Using Gr\"obner Bases}
\noindent
Recall that monomials correspond to points of the lattice $\N^n$.  A
\emph{monomial order} is a total order $\preceq$ of $\N^n$, such that for any
$\alpha,\beta,\gamma\in\N^n$, if $\alpha\preceq\beta$ then
$\alpha+\gamma\preceq\beta+\gamma$.  In other words, a monomial order is a
total order which is compatible with addition of the points in $\N^n$, which
corresponds to multiplication of monomials.

An example of a monomial order is the \emph{lexicographic order}, which is
defined as follows.  First define some ordering on the variables, e.g.,
$x_n\preceq x_{n-1}\preceq\cdots\preceq x_2\preceq x_1$.  Then the
lexicographic order can be defined recursively on the number $n$ of variables:
$x_n^{\alpha_n}\preceq x_N^{\beta_n}$ if and only if $\alpha\leq\beta$, and 
$x_1^{\alpha_1}x_2^{\alpha_2}\cdots x_n^{\alpha_n}\preceq
x_1^{\beta_1}x_2^{\beta_2}\cdots x_n^{\beta_n}$ 
if and only if either $\alpha_1\leq\beta_1$, or $\alpha_1=\beta_1$ and
$x_2^{\alpha_2}\cdots x_n^{\alpha_n}\preceq x_2^{\beta_2}\cdots x_n^{\beta_n}$.

Suppose we are given two polynomials in $n$ variables,
$f_1(x_1,\ldots,x_n)$ and $f_2(x_1,\ldots,x_n)$.  The set of monomials
occurring with nonzero coefficients in $f_i$ is the \emph{support} ${\mathcal
A}_i$ for each $i$.  (Recall that the Newton polytope is the convex hull of the
corresponding lattice points.)  Using the lexicographic order, we can write the
elements of ${\mathcal A}_i$ in a unique way as
$m_{i0},m_{i1},\ldots,m_{i{j_i}}$ such that 
$m_{i0}\prec m_{i1}\prec\cdots\prec m_{i{j_i}}$.  Then we can write the
polynomial equations as
\[f_1(x_1,\ldots,x_n)=a_{j_1}m_{1{j_1}} + \cdots + a_1 m_{11} +
a_0 m_{10}\]
and
\[f_2(x_1,\ldots,x_n)=b_{j_2}m_{2{j_2}} + \cdots + b_1 m_{11} +
b_0 m_{10}.\]
Here $m_{1{j_1}}$ is called the \emph{leading monomial} of $f_1$,
$\alpha_{j_1}m_{1{j_1}}$ is called the \emph{leading term} of $f_1$, and
$\alpha_{j_1}$ is called the \emph{leading coefficient} of $f_1$, and
similarly for $f_2$.  

Now that we have a definite order in which to write the monomials, we can use
long division to divide one polynomial by another.  This is very much like long
division in arithmetic (in fact, in a sense, it's easier, since there's nothing
to guess).  Suppose 
$m_{1{j_1}}=x_1^{\alpha_1}x_2^{\alpha_2}\cdots x_n^{\alpha_n}$ and
$m_{2{j_2}}=x_1^{\beta_1}x_2^{\beta_2}\cdots x_n^{\beta_n}$.
To divide $f_1$ by $f_2$, we would write $f_1$ and underneath it
$a_{j_1}b_{j_2}^{-1}x_1^{\alpha_1-\beta_1}x_2^{\alpha_2-\beta_2}
\cdots x_n^{\alpha_n-\beta_n}f_2$.  When this is not a proper polynomial because
$\beta_j>\alpha_j$ for some $j$, i.e., the monomial $m_{2{j_2}}$ does not
divide $m_{1{j_1}}$, we're already done: the quotient is $0$
and the remainder is $f_2$ itself.  Otherwise, we subtract this from $f_1$ and
write a term 
$a_{j_1}b_{j_2}^{-1}x_1^{\alpha_1-\beta_1}x_2^{\alpha_2-\beta_2}\cdots
x_n^{\alpha_n-\beta_n}$
in the quotient.  The leading term of $f_1$ cancels out, and the leading
monomial of the difference is strictly smaller.  Then we repeat
this process again on this result, adding another term to the quotient,
until either we can't use the leading term of $f_2$ to cancel the leading term
of the result, or the result is zero.  Then this last difference is the
remainder, and we have written down all of the quotient.  The process has to
terminate
because the leading monomials keep getting smaller and smaller.  

Suppose that the two polynomial equations $f_1(x_1,\ldots,x_n)=0$ and
$f_2(x_1,\ldots,x_n)=0$ hold.  From these two polynomial equations we can
derive some more polynomial equations which are logical consequences of them.  
Let $\gamma_i=\max(\alpha_i,\beta_i)$ for $i=1,\ldots,n$, and write
$\gamma=(\gamma_1,\ldots,\gamma_n)$.
Since $a_{j_1}\neq 0$ and $b_{j_1}\neq 0$, we have the following equation:
\[a_{j_1}^{-1}x_1^{\gamma_1-\alpha_1}
x_2^{\gamma_2-\alpha_2}
\cdots
x_n^{\gamma_n-\alpha_n}f_1
-b_{j_2}^{-1}
x_1^{\gamma_1-\beta_1}
x_2^{\gamma_2-\beta_2}
\cdots
x_n^{\gamma_n-\alpha_n}f_2=0.\]
Alternatively, we could use instead the equation
\[b_{j_2}x_1^{\gamma_1-\alpha_1}
x_2^{\gamma_2-\alpha_2}
\cdots
x_n^{\gamma_n-\alpha_n}f_1
-a_{j_1}
x_1^{\gamma_1-\beta_1}
x_2^{\gamma_2-\beta_2}
\cdots
x_n^{\gamma_n-\alpha_n}f_2=0,\]
which does not require that $a_{j_1}\neq 0$ or $b_{j_2}\neq 0$.
We have chosen the polynomials with which to multiply $f_1$ and $f_2$ in order
to cancel the leading terms of $f_1$ and $f_2$.  This polynomial is called
the \emph{$S$-polynomial} $S(f_1,f_2)$ of $f_1$ and $f_2$.  Clearly 
$S(f_1,f_2)(x_1,\ldots,x_n)=0$ also.

Now we can divide $S(f_1,f_2)$ by $f_1$, getting an equation
$S(f_1,f_2)=q_1f_1+r_1$ for polynomials $q_1$ and $r_1$, and then divide $r_1$
by $f_2$, getting $r_1=q_2f_2+r_2$ for polynomials $q_2$ and $f_2$.  (It is
unfortunately the case that the final remainder $r_2$ depends on the order in
which we divided by $f_1$ and $f_2$.)  We have that
$r_2=S(f_1,f_2)-q_1f_1-q_2f_2$, so in particular, $r_2(x_1,\ldots,x_n)=0$.
Thus, if $r_2$ is nonzero, we have a logical consequence of our polynomial
equations, and we can throw it into our polynomial system.  Our polynomial
system is now $\{f_1,f_2,r_2\}$.

Now if we repeat the process, taking $S$-polynomials of pairs of polynomials in
our new system and dividing each $S$-polynomial by all the polynomials in our
new system, we may find more polynomial equations to throw into the system.  It
is a fact from commutative algebra that this process, called \emph{Buchberger's
algorithm}, will always terminate (i.e., finally all the remainders will be
zero), and the (finite) system we have at the end is called a \emph{Gr\"obner
basis}.

Gr\"obner bases have many nice properties, but what will be important to us for
solving polynomial systems is \emph{elimination theory}.  If we compute a
Gr\"obner basis of a polynomial system in the lexicographic order with
$x_n\preceq x_{n-1}\preceq\cdots\preceq x_1$, then those elements of the
Gr\"obner basis involving only $x_n$ will tell us exactly what polynomial
equations in $x_n$ alone which are logical consequences of the polynomial
system.  We can find the roots (if we prefer, only the real roots)
of a polynomial equation in one variable, which
gives us the possible values of $x_n$.  Those elements of the Gr\"obner basis
involving only $x_{n-1}$ and $x_n$ will tell us exactly the polynomial
equations in $x_{n-1}$ and $x_n$ which are the logical consequences of the
polynomial system.  We can substitute in the possible values of $x_n$ we got
before, to get the possible alternative polynomial equations that $x_{n-1}$
alone could satisfy.  Thus we can get the possible values of
$(x_{n-1},x_n)$.  Continuing in this way, we can get all the possible values of
$(x_1,\ldots,x_n)$.  If the system is positive-dimensional, we won't be able
to do this $n$ times.  But the generic finiteness theorem of Harsanyi
\cite{harsanyi} tells us
that a generic game has a finite number of Nash equilibria, i.e., the
associated polynomial systems are zero-dimensional.

Let's use the software package \Singular\ \cite{singular} 
to compute the Gr\"obner basis of a
game of 3 players with 2 pure strategies each.
\begin{verbatim}
~$ Singular
                     SINGULAR                             /
 A Computer Algebra System for Polynomial Computations   /   version 3-0-2
                                                       0<
     by: G.-M. Greuel, G. Pfister, H. Schoenemann        \   July 2006
FB Mathematik der Universitaet, D-67653 Kaiserslautern    \
> ring R=(0,u1100,u1101,u1110,u1111,u2010,u2011,u2110,u2111,
. u3001,u3011,u3101,u3111),(s10,s20,s30,s11,s21,s31),lp;
\end{verbatim}
Here the second parenthesized expression gives the unknowns.  As before we
write {\tt s10} for $\sigma_{10}$ and so forth.
In the first parenthesized expression, the first element denotes the
characteristic of the ring.  To compute over the rational numbers $\Q$ we set
the characteristic to $0$.  The rest of the elements denote \emph{parameters}.
In this case, the parameter {\tt u1ijk} denotes
$u_1(s_{1i},s_{2j},s_{3k})-u_1(s_{10},s_{2j},s_{3k})$, the parameter {\tt
u2ijk} denotes $u_2(s_{1i},s_{2j},s_{3k})-u_2(s_{1i},s_{20},s_{3k})$, and so
forth.  Finally, {\tt lp} means to use the lexicographic order.

Next we specify our polynomial system:
\begin{verbatim}
> poly g1=u1100*s20*s30+u1101*s20*s31+u1110*s21*s30+u1111*s21*s31;
> poly g2=u2010*s10*s30+u2011*s10*s31+u2110*s11*s30+u2111*s11*s31;
> poly g3=u3001*s10*s20+u3011*s10*s21+u3101*s11*s20+u3111*s11*s21;
> poly g4=s10+s11-1;
> poly g5=s20+s21-1;
> poly g6=s30+s31-1;
\end{verbatim}
Finally we ask \Singular\ to compute a Gr\"obner basis:
\begin{verbatim}
> ideal G=g1,g2,g3,g4,g5,g6;
> G = groebner(G);
\end{verbatim}
We ask \Singular\ to display the Gr\"obner basis it computed:
\begin{verbatim}
> G;
G[1]=s30+s31-1
G[2]=s20+s21-1
G[3]=s10+s11-1

G[4]=(u1100*u2011*u3101-u1100*u2011*u3111-u1100*u2111*u3001+
u1100*u2111*u3011-u1101*u2010*u3101+u1101*u2010*u3111+
u1101*u2110*u3001-u1101*u2110*u3011-u1110*u2011*u3101+
u1110*u2011*u3111+u1110*u2111*u3001-u1110*u2111*u3011+
u1111*u2010*u3101-u1111*u2010*u3111-u1111*u2110*u3001+
u1111*u2110*u3011)*s21+
(u1100*u2010*u3111-u1100*u2011*u3111-u1100*u2110*u3011+
u1100*u2111*u3011-u1101*u2010*u3111+u1101*u2011*u3111+
u1101*u2110*u3011-u1101*u2111*u3011-u1110*u2010*u3101+
u1110*u2011*u3101+u1110*u2110*u3001-u1110*u2111*u3001+
u1111*u2010*u3101-u1111*u2011*u3101-u1111*u2110*u3001+
u1111*u2111*u3001)*s31+
(-u1100*u2010*u3111-u1100*u2011*u3101+u1100*u2011*u3111+
u1100*u2110*u3011+u1100*u2111*u3001-u1100*u2111*u3011+
u1101*u2010*u3101-u1101*u2110*u3001+u1110*u2010*u3101-
u1110*u2110*u3001-u1111*u2010*u3101+u1111*u2110*u3001) 

G[5]=(-u1100*u2011*u3011+u1100*u2011*u3111+u1100*u2111*u3011-
u1100*u2111*u3111+u1101*u2010*u3011-u1101*u2010*u3111-
u1101*u2110*u3011+u1101*u2110*u3111+u1110*u2011*u3001-
u1110*u2011*u3101-u1110*u2111*u3001+u1110*u2111*u3101-
u1111*u2010*u3001+u1111*u2010*u3101+u1111*u2110*u3001-
u1111*u2110*u3101)*s11+
(-u1100*u2010*u3001+u1100*u2010*u3011+u1100*u2011*u3001-
u1100*u2011*u3011+u1100*u2110*u3001-u1100*u2110*u3011-
u1100*u2111*u3001+u1100*u2111*u3011+u1101*u2010*u3001-
u1101*u2010*u3011-u1101*u2011*u3001+u1101*u2011*u3011-
u1101*u2110*u3001+u1101*u2110*u3011+u1101*u2111*u3001-
u1101*u2111*u3011+u1110*u2010*u3001-u1110*u2010*u3011-
u1110*u2011*u3001+u1110*u2011*u3011-u1110*u2110*u3001+
u1110*u2110*u3011+u1110*u2111*u3001-u1110*u2111*u3011-
u1111*u2010*u3001+u1111*u2010*u3011+u1111*u2011*u3001-
u1111*u2011*u3011+u1111*u2110*u3001-u1111*u2110*u3011-
u1111*u2111*u3001+u1111*u2111*u3011)*s21*s31+
(u1100*u2010*u3001-u1100*u2010*u3011-u1100*u2011*u3001+
u1100*u2011*u3011-u1100*u2110*u3001+u1100*u2110*u3011+
u1100*u2111*u3001-u1100*u2111*u3011-u1110*u2010*u3001+
u1110*u2010*u3011+u1110*u2011*u3001-u1110*u2011*u3011+
u1110*u2110*u3001-u1110*u2110*u3011-u1110*u2111*u3001+
u1110*u2111*u3011)*s21+
(u1100*u2010*u3001-u1100*u2010*u3011+u1100*u2010*u3111-
u1100*u2011*u3001+u1100*u2011*u3011-u1100*u2011*u3111-
u1100*u2110*u3001+u1100*u2111*u3001-u1101*u2010*u3001+
u1101*u2010*u3011-u1101*u2010*u3111+u1101*u2011*u3001-
u1101*u2011*u3011+u1101*u2011*u3111+u1101*u2110*u3001-
u1101*u2111*u3001-u1110*u2010*u3101+u1110*u2011*u3101+
u1110*u2110*u3001-u1110*u2111*u3001+u1111*u2010*u3101-
u1111*u2011*u3101-u1111*u2110*u3001+u1111*u2111*u3001)*s31+
(-u1100*u2010*u3001+u1100*u2010*u3011-u1100*u2010*u3111+
u1100*u2011*u3001+u1100*u2110*u3001-u1100*u2111*u3001-
u1101*u2010*u3011+u1101*u2010*u3111+u1110*u2010*u3101-
u1110*u2011*u3001-u1110*u2110*u3001+u1110*u2111*u3001+
u1111*u2010*u3001-u1111*u2010*u3101)

G[6]=(-u1100*u2010*u3111+u1100*u2011*u3111+
u1100*u2110*u3011-u1100*u2111*u3011+u1101*u2010*u3111-
u1101*u2011*u3111-u1101*u2110*u3011+u1101*u2111*u3011+
u1110*u2010*u3101-u1110*u2011*u3101-u1110*u2110*u3001+
u1110*u2111*u3001-u1111*u2010*u3101+u1111*u2011*u3101+
u1111*u2110*u3001-u1111*u2111*u3001)*s31^2+
(2*u1100*u2010*u3111-u1100*u2011*u3111-
2*u1100*u2110*u3011+u1100*u2111*u3011-u1101*u2010*u3111+
u1101*u2110*u3011-2*u1110*u2010*u3101+u1110*u2011*u3101+
2*u1110*u2110*u3001-u1110*u2111*u3001+u1111*u2010*u3101-
u1111*u2110*u3001)*s31+
(-u1100*u2010*u3111+u1100*u2110*u3011+u1110*u2010*u3101-
u1110*u2110*u3001)
\end{verbatim}
(We have reformatted the output.)
The Gr\"obner basis has $6$ elements.  The first three elements tell us how to
find {\tt s10} in terms of {\tt s11}, {\tt s20} in terms of {\tt s21}, and {\tt
s30} in terms of {\tt s31}.  The last element {\tt G[6]} is a quadratic
polynomial in {\tt s31}
alone.  We can solve this equation 
to find the possible values of {\tt s31}.  The fourth
element {\tt G[4]} tells us how to obtain {\tt s21} once we have {\tt s31},
and the fifth element {\tt G[5]} tells us how to obtain {\tt s11} once we have
{\tt s31} and {\tt s21}.

Any particular $2\times2\times2$ game is specified by particular values of the
parameters, so we can just substitute them in and solve the resulting system.
Having the Gr\"obner basis gives us important information about how the
geometry of the solution set varies with the parameters.  For instance, if the
coefficient of ${\tt s31}^2$ in {\tt G[6]} vanishes, then this polynomial only
has degree $1$ and hence only one real root.  We can
consider the \emph{discriminant} of the quadratic equation {\tt G[6]}.  Writing
$u^{i}_s=u_i(s_{ij},s_{-i})$, the discriminant becomes:
\begin{eqnarray*}
\biggl(2u^1_{100}u^2_{010}u^3_{111} - u^1_{100}u^2_{011}u^3_{111} -
2u^1_{100}u^2_{110}u^3_{011} + u^1_{100}u^2_{111}u^3_{011}\\
\qquad\qquad- u^1_{101}u^2_{010}u^3_{111} + u^1_{101}u^2_{110}u^3_{011} -
2u^1_{110}u^2_{010}u^3_{101} + u^1_{110}u^2_{011}u^3_{101} +\\
\qquad\qquad2u^1_{110}u^2_{110}u^3_{001} - u^1_{110}u^2_{111}u^3_{001} +
u^1_{111}u^2_{010}u^3_{101} - u^1_{111}u^2_{110}u^3_{001}\biggr)^2\\
- 4\biggl( - u^1_{100}u^2_{010}u^3_{111} + u^1_{100}u^2_{011}u^3_{111} +
 u^1_{100}u^2_{110}u^3_{011} - u^1_{100}u^2_{111}u^3_{011} +\\
\qquad\qquad u^1_{101}u^2_{010}u^3_{111} - u^1_{101}u^2_{011}u^3_{111} -
 u^1_{101}u^2_{110}u^3_{011} + u^1_{101}u^2_{111}u^3_{011} +\\
\qquad u^1_{110}u^2_{010}u^3_{101} - u^1_{110}u^2_{011}u^3_{101} -
 u^1_{110}u^2_{110}u^3_{001} + u^1_{110}u^2_{111}u^3_{001} -\\
\qquad u^1_{111}u^2_{010}u^3_{101} + u^1_{111}u^2_{011}u^3_{101} +
 u^1_{111}u^2_{110}u^3_{001} - u^1_{111}u^2_{111}u^3_{001}\biggr) \\
\biggl( - u^1_{100}u^2_{010}u^3_{111} + u^1_{100}u^2_{110}u^3_{011} +
u^1_{110}u^2_{010}u^3_{101} - u^1_{110}u^2_{110}u^3_{001}\biggr).
\end{eqnarray*}
The set of payoff functions where the discriminant is zero is a real algebraic
variety in the space of $2\times2\times2$ games, the \emph{discrimant variety}.
It partitions the space of $2\times2\times2$ games into a region where the
discriminant is positive, in which case the polynomial system has two real
roots, and a region where the discriminant is negative, in which case the
polynomial system has no real roots.

The reader was already familiar with the discriminant of the quadratic formula,
but the same phenomenon will happen with equations of higher degrees in more
variables \cite{GKZ}.  In this case there could be several discriminant
varieties, with various implications about the geometry of the solution set.
The {\sf SALSA} team of INRIA Rocquencourt and LIP6 in France has produced a
Maple package {\sf DV} which will analyze the discriminant varieties for a
parametric polynomial system \cite{LR06}, and Antonio Montes \cite{montes} has
independently produced another such Maple package {\sf DisPGB}.  (These
analyze the implications of the discriminantal equations being zero or nonzero
for the complex solutions of the polynomial system.  In the example above, we
considered in addition the implications of the discriminant being positive or
negative on the real solutions.)

Specifically, Gabriela Jeronimo, Daniel Perrucci, and Juan Sabia have recently
explained how to obtain a parametric representation of the totally mixed Nash
equilibria.  (As we know, this means we can get a parametric representation of
all the Nash equilibria by considering various possible supports.)  They give
polynomial-time algorithms for describing the set of totally mixed Nash
equilibria, using multihomogeneous resultants.  The \emph{resultant} of a
polynomial system is a polynomial equation in the coefficients of the system
which must hold in order for the system to have a root.  Resultants are a key
tool in the solution of polynomial systems.  The parametric equations
characterize the geometry of the space of games and provide an effective 
method for finding all the roots.  We look forward eagerly to the
implementation of these algorithms in a software package, which we hope will
lead to many new insights in game theory.

\section{Finding All Nash Equilibria in \Gambit}
\noindent
The \Gambit\ software package incorporates a variety of tools for finding Nash
equilibria and studying other properties of games.  We will discuss here the
version of Jaunary 6th, 2006, which was the latest released version at the time
of this writing.  

The \Gambit\ source includes a procedure to call the {\sf
Pelican}
software for polyhedral homotopy continuation which was written by Birk Huber
\cite{hubersturmfels}.  However, as of this release, although \Gambit\ does use
homotopy continuation for computing the logistic Quantal Response Equilibrium
correspondence, as well as the tracing procedures mentioned above, it does not
use polyhedral homotopy continuation to solve the polynomial systems in order
to enumerate all Nash equilibria.  The \Gambit\ source also includes code to
solve the polynomial systems by Gr\"obner basis techniques, but \Gambit\ does
not actually use these either.  

{\sf Pelican}, written in 1995, is no longer actively maintained, so it's not
surprising that \Gambit\ doesn't use it.  There are now a few other packages
for polyhedral homotopy continuation, however, and Jan Verschelde continues to
actively develop \PHC{pack} in particular.  \PHC{pack} also includes a C
interface.  Furthermore, the choice of a
factorizable starting system with manifest roots which we have described here
avoids the more delicate and computationally complex aspects of polyhedral
homotopy continuation, leading to an efficient solution method in practice.

The \Gambit\ Gr\"obner basis code, which was written specifically for \Gambit,
was abandoned due to ``numerical instability''.  Using Gr\"obner bases for
numerical root finding is indeed a delicate procedure, since computing a
Gr\"obner basis at machine precision involves many intermediate multiplication
and division operations, which makes it difficult to keep errors from
accumulating beyond a tolerable level.  Of course, at the end one also has to
find the roots of a univariate polynomial, but many carefully written numerical
libraries offer routines to do this.  The bigger problem is that after all the
preceding computation to arrive at the Gr\"obner basis, the coefficients of the
univariate polynomial we're solving may be wrong, giving us a root which is
wrong; the coefficients of the bivariate polynomial into which we substitute
the root of the univariate polynomial may also be wrong, giving us a root which
is even more wrong; and so forth.  It is for these reasons that we too
recommend Gr\"obner bases for geometric insight into how the structure of the
set of Nash equilibria may vary over the space of all games of a particular
format, rather than as a practical tool for finding all the Nash equilibria of
one particular game.  By contrast, homotopy continuation is relatively more
mature and well-understood as a numerical technique.

Instead, \Gambit\ solves the polynomial systems defining the Nash equilibria by
subdividing the product of simplices (where each simplex is the subset of
$\R^{d_i}$ defined by $\sigma_{i1}+\cdots+\sigma_{i{d_i}}=1$, $\sigma_{ij}\geq
0$ for all $j$) into small boxes and looking at the Taylor series of each
polynomial in the system.  We can, for instance, evaluate the $D$ polynomials
$(f_1,\ldots,f_D)$ at one corner $\sigma$ of such a box, giving a vector
$(y_1,\ldots,y_D)$.  Then any point in the box is no further than the
diagonally opposite corner, and we can plug this distance into the Taylor
series of the polynomials about the corner to determine a bound on how far the
image of any other point $\sigma'$ of the box can lie from $(y_1,\ldots,y_D)$.
If this bound is less than $||(y_1,\ldots,y_D){||}_2$, then no point in the box
can be a root of the polynomial system.  If it is not, then \Gambit\ looks for
a root in the box using Newton's method.  If such a root is found, then
\Gambit\ tries to see whether it can determine that no other roots may exist
within the box, again using the Taylor series of the polynomials.  If it can,
then it is done with this box, but if not, then it subdivides the box into
$2^D$ smaller boxes and looks at those.

This method has the advantage that it will only find real roots which satisfy
the nonnegativity constraints.  However, it does not scale well with higher
dimensions.  Using polyhedral homotopy continuation with factorizable
starting systems, finding all Nash equilibria of games of much larger formats
should become practical.  However, it may be useful to use a similar Taylor
series technique on the augmented system (the one including the variable $t$)
along the way to each root of the target system, to see
whether we can derive a bound ensuring that the target root will not be real
and nonnegative (i.e., a Nash equilibrium).  In that case we can abandon this
particular path.  Here we are just travelling along a $1$-dimensional interval
$[0,1]$, so the problem of having to subdivide into an exponential number $2^D$
of smaller boxes does not arise.

In this release \Gambit\ provide an option to compute Nash equilibria via
heuristic search on the supports of the game \cite{PorterNudelmanShoham},
through an implementation contributed by Litao Wei.  Heuristic search is
complementary to using polynomial algebra to find Nash equilibria.  It can tell
us which of the many possible supports to look at first.  The choice of support
tells us which polynomial systems to solve (possibly none, if we find a pure
Nash equilibrium on the support and are satisfied with not looking any further).

\section{Conclusion}
In this paper we have described the polynomial systems which characterize the
Nash equilibria of a game.  We have explained how to construct and solve a
factorizable start system and then use polyhedral homotopy continuation to
solve games of a given format.  We have also explained how to use Gr\"obner
bases to gain insight into how the geometry of the solution set of these
polynomial systems varies over the space of games of a particular format.
Finally, we have reviewed the current use of \Gambit\ for finding all Nash
equilibria of a game.  We suggest that \Gambit\ may be able to find all Nash
equilibria of games of larger formats than is currently possible, by
incorporating polyhedral homotopy continuation from factorizable nondegenerate
start systems.  Alternatively, or in addition, \Gambit\ may implement the
algorithm of \JeronimoPerrucciSabia\ to find all Nash equilibria of a
game.  We sincerely hope that the possibility of analyzing larger games will
enable game theorists to make more realistic models of strategic interaction.

\section{Acknowledgements}

Our earlier paper \cite{dattaISSAC} contains much of the material which is
surveyed more expansively here.  We would like to express our gratitude to the
following for generously taking the time to personally discuss with us the use
of their software packages: Andrew McLennan and Ted Turocy ({\sf Gambit}
\cite{gambit}), Gert-Martin Greuel ({\sf Singular} \cite{singular}), and Jan
Verschelde ({\sf PHC} \cite{phc}).  We would also like to thank Gabriela
Jeronimo for sending us a preprint of her paper with Daniel Perrucci and Juan
Sabia, and Andrew
McLennan for suggesting she do so.  We would like to
thank Richard Fateman and Bernd Sturmfels for supervising the research leading
up to that paper, during which the author was partially supported by NSF grant
DMS 0138323.  We would also especially like to acknowledge our debt to Bernd
Sturmfels, especially for teaching us about the application of polynomial
algebra to Nash equilibria, in the lectures leading to the book
\cite{sturmfels}.

\bibliographystyle{abbrv}
\bibliography{JETSurvey,computenash,thesis}

\end{document}